\documentclass[12pt,a4paper,draft]{article}
\usepackage{russian}
\unitlength=1.00mm\linethickness{0.4pt}
\newfont{\blef}{cmff10}
\renewcommand{\S}{{\blef\begin{picture}(2,0)\put(0,-0.27){S}\put(0.08,0.73){S}\end{picture}}}

\newcommand{\st}{\scriptstyle}
\newcommand{\sx}[2]{\stackrel{[#1]}{#2}}
\newcommand{\sy}[1]{\stackrel{\circ}{#1}}
\newcommand{\G}{{\rm G}} \newcommand{\PL}{\prod\limits}
\newcommand{\Le}{\left[\begin{array}{l}} \newcommand{\R}{\end{array}\right]}
\newcommand{\dt}[1]{\hbox to #1{{}\dotfill {}}}
\newcommand{\q}{q^{}} \newcommand{\pe}{p^{}} \newcommand{\er}{r^{}}
\newcommand{\x}{x^{}} \newcommand{\y}{y^{}}
\newcommand{\tri}{\makebox(0,0)[cc]{$\st\triangle$}}
\newcommand{\ti}{\makebox(0,0)[cc]{$\times$}}
\newcommand{\frm}[1]{\framebox(1.65,1.65)[cc]{#1}}
\newcommand{\fra}{\rule{1.70mm}{1.70mm}}
\newcommand{\ddt}[1]{\begin{picture}(20,3)
\multiput(0,0)(5,0){#1}{$-\ $}\end{picture}}
\newcommand{\abzats}{\hangindent=\parindent\hangafter=-1\noindent}

\begin{document}
\thispagestyle{empty}

\begin{center} JOINT INSTITUTE FOR NUCLEAR RESEARCH \end{center}
\begin{center} Laboratory of Computing Techniques and Automation
\end{center}

\vspace*{4cm}

\begin{center} G.A.Emel'yanenko, M.G.Emelianenko\footnote[0]{Computational
mathematics and cybernetics department, Moscow State University.},
T.T.Rakhmonov, E.B.Dushanov, \\ G.Yu.Konovalova\footnote{International
University of Nature, Society, and Human Being ``Dubna"}\end{center}

\vspace*{1cm}

\begin{center}{\Large On efficiency of}\end{center}
\begin{center}{\Large critical-component method for solving}\end{center}
\begin{center}{\Large singular and ill-posed systems of}\end{center}
\begin{center}{\Large linear algebraic equations}\end{center}

\vspace*{2cm}

\begin{center} Abstract \end{center}

Results are expounded for the investigation of efficiency of the
critical-component direct method for solving {\it degenerate} and
{\it ill-posed} systems of linear algebraic equations.

\vspace*{5cm}

\begin{center} DUBNA,1998  \end{center}

\newpage
\subsubsection{Introduction}
\abzats
In this paper, we present results of studies of efficiency of the critical-
component direct method proposed in [$1\div3$] for solving {\it degenerate}
and {\it ill-posed} systems of linear algebraic equations
			$$AZ=F,\leqno(1.1)$$
where $A$ is a square matrix of the general form with real elements
$a^{}_{ij}, A=\{a^{}_{ij}\}, Z$ is an unknown vector with coordinates
$z^{}_j, Z=\{z^{}_j\}$, and $F$ is a known vector with coordinates $f^{}_i,
F=\{f^{}_i\}, i,j=1,2,...,m.$ It is shown that for systems like (1.1) the
critical-component method makes it possible to numerically determine the only
normal pseudosolution $(Z^+\!\!=A^+F)\!:\
||AZ^+\!\!-F||=\!\!\!\inf\limits_{Z\in Z^{}_{\st A}}\!\!\!||AZ-F||,
||Z^+\!||=\!\!\!\inf\limits_{Z\in Z^{}_{\st A}}\!\!\!||Z||$, where $Z^{}_{\st
A}$ is a set of all pseudosolutions to system (1.1), and to obtain the unique
matrix $A^+$, pseudoinverse of $A\!:\ ||A^+\!A-E||=\!\!\!\inf\limits_{
\sy{A}\!{}^{-1}\in\Omega^{}_{\st A}}\!\!\!||\sy{A}\!{}^{-1}A-
E||, ||A^+\!||=\!\!\!\inf\limits_{\sy{A}\!{}^{-1}\in\Omega^{}_{\st A}}
\!\!\!||\sy{A}\!\!{}^{-1}||, A^+\!A=AA^+$, where $E$ is a unit matrix and
$\Omega^{}_{\st A}$ is a set of all $\sy{A}\!{}^{-1},$
pseudoinverse of $A$. In this case, even if the problem (1.1) is
substantially ill- posed, the quantities $Z^+$ and $A^+$ are stable to small
changes of input data $(A,F)$. Comparative analysis of results of the
numerical solution performed for a large number of problems like (1.1) both
by the new method and by those known earlier shows that the critical-component
method is on the average more effective than any method compared to it.
When $\det A\ne0$ and a system is well-posed, the normal pseudosolution
$Z^+$ of system (1.1) derived by the critical-component method coincides
with its usual solution $Z$, and $A^+=A^{-1}$ is a matrix inverse of $A$.
One of the main problems in numerical solution of ill-posed systems of
algebraic equations is well-known [4,5,6]: there can be large changes in
the solution, beyond the scope of admissible values, corresponding to small
changes in the matrix of a system or/and its right-hand side. The above
breakdown of continuity of the inverse mapping
$Z=A^{-1}F$, if $A^{-1}$ exists, is caused by a great norm
$||A^{-1}||$ and, as a result, by large $\mu=$ cond $A$,
the condition number of  the system matrix
($\mu=||A||\cdot||A^{-1}||$, if $\det A\ne0$ and $\mu=\infty$, if $\det
A=0$, where $||\cdot||$ are the corresponding norms), i.e. even for an
exactly given vector $F$ a negligible relative error in calculating
$A^{-1}$ can produce a large distortion of the searched vector $Z$.
This effect is to be taken into account since realistic calculations are
carried out with a certain finite accuracy and , besides, sometimes one
knows not the exact system $AZ=F$, but only a system
$\tilde A\bar Z=\tilde F$, approximate of it, which obeys the inequalities
$||\tilde A-A||\leq h^*$ and $||\tilde F-F|| \leq\delta^*$ (the meaning of
norms is defined by the character of a problem). The numbers $h^*>0$ and
$\delta^*>0$, specifying the norms of deviations of approximate data $(\tilde
A, \tilde F)$ from the exact ones $(A,F)$ of problem (1.1) $(h^*\leq
h_0+h_1,\delta^*\leq\delta_0+\delta_1,h_0\geq0,h_1>0,\delta_0\geq0,\delta_1>0)$,
are sums of $(h_0,\delta_0)$,proper model (complete) errors of problem
(1.1) and of $(h_1,\delta_1)$, round-off errors [7,8] when writing the
data into the computer memory. Since there are, thus, infinitely many
systems (1.1) with the input data $(A,F)$, indistinguishable within the
accuracy $(h^*,\delta^*)$, we can speak only about deriving an approximate
solution to system (1.1). As a result, difficulties may arise in numerical
computations for some systems of equations (1.1) with square matrices when
answering the following questions:
\hspace*{2cm}
\begin{itemize}
\item[]-- is the system degenerate ``within accuracy $(h^*,\delta^*)$"
ill-posed?\footnote[0]{Systems degenerate ``within accuracy
$(h^*,\delta^*)$" are not always ill-posed (as example is system
(1.1) with $A=A^T$, singular (eigen)values
$\mu_1=\mu_2=...=\mu_m=10^ {-6},$ determinant $\det A=10^{-6m}$ and the
condition number cond $A=\mu=||A|| \cdot||A^{-1}||=1$).} and\\
\\
-- is a given system ill-posed by virtue of its being degenerate or is it
nondegenerate but ill-posed?\end{itemize}
\hspace*{1cm}

Indeed, if the system $AZ=F$ with a square matrix is degenerate, then $\det
A=0$, i.e., the matrix $A$ has some of its eigenvalues equal to zero.
But if  $\det A\ne0$, and the system is ill-posed, then the normal matrix
$A^TA$ has some eigenvalues only close to zero $\mu^2_1,...,\mu^2_m$
($|\mu_i|$ are singular values of the matrix $A$).
\noindent Consequently, systems of linear algebraic equations with square matrices,
which are ill-posed and degenerate "within a given accuracy $(h^*,\delta^*)$"
may turn out to be indistinguishable in the process of computations.
Besides, the problems (1.1) and $\tilde A\bar Z=\tilde F$ can be inconsistent
if one defines the criterion of consistency [9] determined by accuracies
$(h^*,\delta^*)$. It may also happen that $\det A=0$ (or $\det\tilde
A=0$), i.e. system (1.1) (or $\tilde A\bar Z=\tilde F$) has an infinite
number of solutions. Then, there arises the question: what is to be
understood by the numerical solution to the initial system $AZ=F$.
There are various conceptual approaches to solve this problem (see, for
instance, reviews given in [4,6,10], etc.).\\

If one takes advantage of the regularization [4], the solution $Z^+$ to the
system $AZ=F$ (1.1) will be the regularized normal pseudosolution
$Z^{\alpha}$ that minimizes the discrepancy $||\tilde A\bar Z-\tilde F||$ on
the set of all its pseudosolutions $Z^{}_{\st A}$ if
$||Z^{\alpha}||=\inf\limits_{\bar Z \in Z^{}_{\st A}}||\bar Z||$ and
$Z^{\alpha}$ is stable to small variations in $(h^*,\delta^*)$ of
input data $(A,F)$. The parametric vector $Z^{\alpha}$ is directly computed
by solving the sequence of normal systems of equations $(\tilde
A{}^T\tilde A+\alpha E)Z^{\alpha}=\tilde A{}^T\tilde F$ with the aim of a
more accurate iterative determination of the minimum of quadratic functional
$M^{\alpha}[\bar Z,\tilde F,\tilde A]=||\tilde A\bar Z-\tilde F||^2+\alpha
||\bar Z||^2$ with the regularization parameter $\alpha(\alpha>0)$,
determined from the discrepancy, i.e., from the condition $||\tilde
AZ^{\alpha}-\tilde F||=\delta_*$, where $\delta_*(\delta_*>0)$ is a
numerical function of $(h^*,\delta^*)$ and $Z^{\alpha}$ [4,5,6].\\

The other group of numerical methods of solving the problem (1.1) rely on
searching for the generalized matrix $A^+$, which is (pseudo)inverse of
$\tilde A$, either by the method of singular decomposition ($\tilde A=U\Sigma
V,$ where $U$ and $V$ are orthogonal matrices, $\Sigma$ is a diagonal matrix,
whose elements are singular numbers
$|\mu_1|\geq|\mu_2|\geq\cdots\geq|\mu_m|\geq0$ of the matrix $A$, and
$A^+=V^T\Sigma^+U^T$), or by some
other method [7,9,10,11]. Common to both of the approaches is that in their
program realization they solve (each by its own means and with its own
efficiency) the problems of minimization of norms $||\tilde A\bar Z-\tilde
F||$ and $||\bar Z||$ and of the continuous dependence of the solution $Z^+$
on small changes in $(h^*,\delta^*)$ of input data $(A,F)$. Here it is set
that $\mu=$ cond $A=||\tilde A||\cdot||A^+||$, and the main problem now is a
stable calculation of the rank of  $\tilde A$ [7,9].\newpage

Conceptually, the critical-component method can be attributed to the second
of indicated groups of methods. It is based on the idea of constructive
search (under the condition that matrix and vector norms are consistent:
$||Z^+||\leq||A^+||\cdot||\tilde F||$;  and  the matrix norm is induced
by the vector norm:
$||A^+||=\sup\limits_{||\tilde F||\ne0}(||A^+\tilde F||/||\tilde F||)$ [5,9])
of an optimal representation for the matrix $A^+$, pseudoinverse of the
matrix $\tilde A$, in the process of decomposition of system (1.1) into
subsystems, whose solution is stable to errors \footnote[0]{Throughout we use
the notation:  $\varepsilon_1 (\varepsilon_1>0)$ is the modulus of relative
error of the arithmetic of computer operations with real numbers with a
floating point; $\varepsilon_0(\varepsilon_0>0)$ is the modulus of absolute
error of the computer zero $\theta$, i.e. of any small real number (except
for 0) from the interval $\theta\in(0-\varepsilon_0,\varepsilon_0+0)$, where
0 is the usual zeroth element of the real axis. If
$\theta\in(0-\varepsilon_0, \varepsilon_0+0)$ and $\theta\ne0$, it is
accepted that $\theta=0$ [8,9].  Using constants $\varepsilon_1$ and
$\varepsilon_0$, one can estimate [7] errors of arrangement (writing) of the
real $[m,m]$ matrix $A$ and $m$-dimensional vector $F$ in the computer memory
in the form $||A_{comp}-A||_E\!\leq\!(\varepsilon_1
||A||_E+\varepsilon_0m\equiv h_1),
||F_{comp}-F||_E\!\leq\!(\varepsilon_1||F||_
E+\varepsilon_0\sqrt{m}\equiv\delta_1)$, where $||\cdot||_E$ are the
Euclidean norms of matrices and vectors, and $h_1>0,\delta_1>0$.}
$(\varepsilon_1,\varepsilon_0)$ and small $(h^*,\delta^*)$ changes of input
data $(A,F)$. High efficiency of the critical-component method is
provided by its basic constituents:
\begin{itemize} \item[]-- the reduction, stable to errors
$(h^*,\delta^*;\varepsilon_1,\varepsilon_0)$, of system (1.1) to
two-(tri)diagonal systems; \\
-- generalized processes $\{\Lambda,\G\}$, stable to errors
$(\varepsilon_1,\varepsilon_ 0)$ [14], for
calculating ratios of upper (lower) corner minors of triangular matrices
which allow one, accurate within constants ($\varepsilon_1$ and
$\varepsilon_0$) of the computer arithmetic, to determine the structure and
diagonal elements of matrices that are inverse of them (introduced in [12,13]);  \\
-- the algorithm of optimal (with $(\varepsilon_1,\varepsilon_0)$)
decomposition of the system $\tilde A\bar Z=\tilde F$ into well-posed
subsystems; \\  -- the algorithm of optimal sewing of the solution
$Z^+$ to the system $\tilde A\bar Z= \tilde F$ from well-posed
subspace solutions.  \end{itemize}

In what follows, along with problem (1.1) of the general form, we will
consider the problems of numerical solution of degenerate and ill-posed
systems of linear algebraic equations
\begin{flushleft}(1.2)\hspace*{7cm}$C_3X=Y,$ \\
(1.3)\hspace*{7cm}$C_2\hat X=\hat Y$ \end{flushleft}
with square real matrices $C_3$ and $C_2$ of order $m$,
of the tridiagonal and two-diagonal form respectively:
$$C_3=\Le\q_1\ \er_2\\ \pe_2\ \q_2\ \er_3\\ \ \ \ddots\ddots\ddots\\ \quad\
\pe_{m-1}\q_{m-1}\er_m\\ \hspace{4em}\ \pe_m\ \q_m\R,\quad C_2=\Le\q_1\
\er_2\\ \quad\ \q_2\ \er_3\\ \qquad\ddots\ddots\\ \qquad\ \ \q_{m-1}\;\er_m\\
\hspace{5em}\q_m\R,\leqno (1.4)$$
where $X=(x_1, x_2, ..., x_m)^T$ and $\hat X=(\hat x_1, \hat x_2, ..., \hat
x_m)^T$ are unknown vectors, and $Y=(y_1, y_2, ..., y_m)^T$ and
$\hat Y=(\hat y_1, \hat y_2, ..., \hat y_m)^T$ are given $m$-dimensional
vectors, $\{q_i\}_{i=1}^m$ are diagonal elements and $\{p_i,r_i\}_{i=2}^m$
are sub(off)diagonal elements of matrices $C_3$ and $C_2$, respectively.
Without loss of generality, we consider systems (1.3) only with the right
two-diagonal matrix.

Since these problems are a particular case of problem (1.1), all said above
applies also to problems (1.2) and (1.3), whose solutions
$X^+$ and $\hat X\!{}^+$ are constructed more easily than $Z^+$. Therefore,
in the course of program realization of the above conceptions of solution of
problem (1.1), the initial stage [4,5,7,9,13,15] consists in its reduction
to problems (1.2) and (1.3), i.e.
$$\left\{\!\!\begin{array}{l}C_2(Q^T\!Z)\!=\!PF,\;\mbox{where}\;C_2\!=PAQ\;
\mbox{is a two-diagonal matrix, if}\;A\ne A^T,\\ C_3(Q^T\!Z)\!=\!Q^TF,\;
\mbox{where}\;C_3\!=Q^T\!AQ\;\mbox{is a tridiagonal matrix, if
}\;A=A^T.  \end{array}\!\!\right.\leqno(1.5)$$
Here $U^TU=E=UU^T, U\!: Q,P$ are matrices of reflections or rotations. The orthogonal
transformations (1.5) stable\footnote[0]{It is known [5,7,9] that the
Euclidean and spectral norms of matrices are invariant (theoretically) under
the orthogonal transformations (1.5), i.e. there hold the equalities:
$||C_2\!||_E\!=||PAQ||_E\!=||A||_E,||C_3||_E\!=||Q^TAQ||_E\!=||A||_E;$
$||C_2\!||_2\!=||PAQ||_2\!=||A||_2,||C_3\!||_2\!=||Q^TAQ||_2\!=||A||_2$ and
$||P||_2\!=||Q||_ 2\!=1,||P||_E\!=||Q||_E\!=\sqrt{m}$. As a result,
$\mu=$ cond $A=$ cond $C_3$ (or cond $C_2$). Here $||PF||_E\!=||\hat
Y||_E,||Q^TF||_E\!=||Y||_E; ||A||_2$ is a norm induced by the Euclidean
vector norms $||Z||_E$ and $||F||_E$; $M(A)=m\!\cdot\!\!\max\limits_{1\leq
i,j\leq m}\!\!\!|a_{ij}|$ and $||A||_E$ are norms consistent with norms
$||Z||_E$ and $||F||_E$. However, in real computations in process (1.5)
of the reduction of system (1.1) to form (1.2) or (1.3), using the
Houscholder $U$ transformations (reflections), we obtain the estimates
[7] \begin{flushleft}
$||(C_2)_{comp}\!-C_2\!||_E=||(PAQ)_{comp}\!-PAQ||_E\leq\left\{\!
\left[\!\left(\!\frac{(2m-3)\varepsilon_r}{1-(m-2)\varepsilon_r}||A||_E+
\frac{(2m-3)\sqrt{m}\cdot0_r}{1-(m-2)\varepsilon_r}\right)\!\equiv\!f_2(m)
||A||_E\!\right]\!\equiv\!h_2\!\right\}\!,$\\ $||(PF)_{comp}-PF||_E\!\leq(\varepsilon_r
||F||_E+0_r)\equiv\delta_2$ and\\ $||(C_3)_{comp}-C_3\!||_E\!=||(Q^TAQ)_{comp}-
Q^TAQ||\!\leq\!\left\{\!\left[\!\left(\!\frac{(2m-4)\varepsilon_r}{1-
(m-2,5)\varepsilon_r}||A||_E+\frac{(2m-4)\sqrt{m}\cdot0_r}{1-(m-2,5)\varepsilon_r}
\right)\!\equiv\!f_3(m)||A||_E\!\right]\!\equiv\!h_2\!\right\}\!,$\\
$||(Q^TF)_{comp}-Q^TF||_E\leq(\varepsilon_r||F||_E+0_r)\equiv\delta_2$, where
$\varepsilon_r\sim29\varepsilon_1$ and $0_r\sim(2m+2\sqrt{m})\varepsilon_0,
h_2>0,\delta_2>0$.\end{flushleft}
Similar inequalities could also be written for $\tilde A, \tilde F,$ where
matrix $\tilde A$ and vector $\tilde F$ differ from $A$ and $F$ by
simultaneous inclusion of inherited errors and errors of writing into the
computer memory. From the above inequalities it follows that problems
$AZ=F$ and $\tilde A\bar Z=\tilde F$ are continuous with respect to the
orthogonal transformations (1.5). Though the inherited errors
$(h_0,\delta_0)$, if known, are, as a rule, much larger than the total
$(h_1+h_2,\delta_1+\delta_2)$ effect of the errors of writing
and transformations (1.5), the latter can influence the character (degree)
of problem (1.2) or (1.3) being well-/ill-posed. The cited monographs contain
also simplified estimates for errors $h_2$ and $\delta_2.$} to errors
$(h^*,\delta^*;\varepsilon_1,\varepsilon_0)$, do not often improve the nature
of the problem being ill-/well- posed. Ill-posed systems of type
(1.1) sometimes can numerically be reduced to ill-posed systems of type (1.2)
and (1.3), with the notation in (1.5): $X=Q^TZ,Y=Q^TF$ and $\hat X=Q^TZ,\hat
Y=PF$.  Therefore the basic problem is numerical solution of such degenerate
and ill-posed systems. Once the vectors $X$ and  $\hat X$ are obtained, we
determine the solution to system (1.1), vector $Z$, in the form
$$Z=QX\mbox{ and  }Z=Q\hat X.\leqno(1.6)$$
Numerical solution of ill-posed systems
(1.2),(1.3) with tridiagonal and upper two-diagonal matrices can be best
realized by the following methods [4,7,9] : the inverse substitution with
normalization, regularization, a singular decomposition with exhaustion.
In sect.3, we present (in particular) the results of comparison between
computations performed by these methodes and by the new one.\newpage

\subsubsection{Critical-component method for numerical solution of degenerate
and ill-posed systems of linear algebraic equations with tri- and
two-diagonal matrices}

\abzats
Below we formulate the theorem according to which one can numerically obtain
the only stable non-iterated normal pseudosolution $X^+$ of the system of
linear algebraic equations of the general form (1.2), stable to errors
$(\varepsilon_1,\varepsilon_0)$ and $(h,\delta)$, by the critical-component
method.

The vector $X^+$ and the representation for the matrix consistent with it
$(C^+_3\equiv B)$, pseudoinverse to $C_3$, are determined as functions of
stably computed vector $\sy{X}$  (a regular component of $X^+$) and
matrix $\sy{B}$ (a regular component of $C^+_3$). In contradistinction to the
problem of computation of singular numbers of matrices $C_3$ being unstable
in nature, the critical-component method is stable owing to the stable
processes of computation
of the ratios of upper (lower) corner minors $\{\Lambda,\G\}$ of this
matrix. Thus, the method of solution based on the search for a non-parametric
stable component of the pseudoinverse matrix [7,9] found one more argument
for its being efficient (contrary to conclusions of perturbation theory
according to which $X^+$ and $C^+_3$ are not valid for computer calculations).

{\bf Theorem.} Let $C_3X=Y$ be either a degenerate or an ill-posed system of
linear algebraic equations with a square, of order $m$, real tridiagonal
matrix of the general form $C_3$ (1.4). Also, let the system
$\tilde C_3\bar X=\tilde Y,$ where $||\tilde C_3-C_3||\leq h$ and $||\tilde
Y-Y||\leq\delta$, being an image of the system $C_3X=Y$ in the computer
memory, be ill-posed but nondegenerate. Then the only pseudosolution $X$ of
the system $C_3X=Y$ that is minimal in norm $(||X^+||=\min)$, obeys the
condition of the norm of discrepancy being minimal
$(||\tilde C_3X^+-\tilde Y||=\min)$, and is stable to computation errors
$(\varepsilon_1,\varepsilon_0)$ and to small changes $(h,\delta)$ of the
input data $(C_3,Y),$ can numerically be obtained by the following direct
critical-component method\footnote[0]{Here $h\leq h_0+h_1+h_2$ and
$\delta\leq \delta_0+\delta_1+\delta_2$ if the system $\tilde C_3\bar
X=\tilde Y$ is a reduced image of the system $AZ=F$; and $h\leq\bar
h_0+\bar h_1, \delta\leq\bar\delta_0+\bar\delta_1,$ where $(\bar
h_0\geq0,\bar\delta_0\geq0)$ are hereditary errors and $(\bar
h_1>0,\bar\delta_1>0)$ are errors of writing the  system $C_3X=Y$ into
the computer memory if system (1.1) is initially of form (1.2).

Since numerical solution is derived for the system $\tilde C_3\bar X=\tilde
Y$ that is, within accuracy $(h,\delta)$, indistinguishable from the
system $C_3X=Y$, for simplicity of the notation, the very algorithm of
numerical method and its proof are given in the notation of the system
$C_3X=Y$, i.e., without ``$\sim$", if this does not cause misunderstanding.
The requirement $\det\tilde C_3\ne0$ of the theorem will be removed later.
$X^+=(x^+_1,x^+_2,...,x^+_m)^T$.}:

\begin{flushleft}
\hspace*{13mm}{\it Start of computations:}\\ \hspace*{13mm}$k=1, i=m;$\\
(2.1)\hspace*{5mm}$l_k=i;$\\
(2.2)\hspace*{5mm}$\sx{k}{x}^{}_i=\sum\limits^{l_k}_{\xi=1}\sx{k}{B}^{}_{i\xi}\y_{\xi},\
\sx{k}{\phi}^{}_i=\left\{\begin{array}{l}0,\mbox{ if
}k=1,\\-\sx{k}{B}^{}_{il_ {\st k}}\er_{l_{\st k}+1}x^+_{l_{\st k}+1},\mbox{
if }k>1;\end{array}\right.$\\ \hspace*{13mm}if $i=l_k$, then (2.5), otherwise
(2.3);\\ (2.3)\hspace*{5mm}if $|\sx{k}{\phi}^{}_i|<1/\varepsilon_1$, then
(2.4), otherwise $k=k+1$ and (2.1);\\ (2.4)\hspace*{5mm}$j=i+1,
\sx{k}{x}^{}_{l_{\st k}+1}=0;$\\
\hspace*{13mm}$\Phi^{}_j=\left\{\begin{array}{l}|\y_j|-|\pe_j\sx{k}{x}^{}_{j-1}+\q_j
\sx{k}{x}^{}_j+\er_{j+1}\sx{k}{x}^{}_{j+1}|,\mbox{ at }|\y_j|\leq 1,\\
1-|\pe_j
\sx{k}{x}^{}_{j-1}+\q_j\sx{k}{x}^{}_j+\er_{j+1}\sx{k}{x}^{}_{j+1}|/|\y_j|,
\mbox{ at }|\y_j|>1;\end{array}\right.$\\
\hspace*{13mm}if $|\Phi^{}_j|\leq 2\varepsilon_1$, then (2.5), otherwise
$k=k+1$ and (2.1);\\
(2.5)\hspace*{5mm}$x^+_i=\sx{k}{x}^{}_i+\sx{k}{\phi}^{}_i;$\\
\hspace*{13mm}if $i=1$, computations are over, otherwise $i=i-1$ and (2.2);\\
\hspace*{13mm}{\it End of computations.} \end{flushleft}

Here:

$\sx{k}{B}^{}_{ij}$ ($l_{k+1}\leq i\leq l_k,\ 1\leq j\leq l_k$ and
$k=1,2,...,n$) are elements of submatrices $\sx{k}{B}$ of the matrix
$\sy{B}=\sy{C}\!{}^ {-1}_3$ that is inverse of a well-posed matrix
$\sy{C}_3$ of the form:  $$\begin{array}{r}(2.6)\\ \\ \sy{C}_3=\\
\end{array}\!  \!\!\!\Le\!\!\!\Le\!\!\!\q_1\ \er_2\\ \!\!\!\pe_2\ \q_2\
\er_3\\ \quad\dt{1cm}\\ \hspace{7mm}\pe_{l_{\st n}}\q_{l_{\st n}}\!\!\!\R\!\!
\begin{array}{l}\\=\sx{n}{C}\!{}^{(l_{\st n+1}+1=1)}_{l_{\st n}}\\ \\0\end{array}\\
\hspace{13mm}\ddt{5}\\ \hspace{1cm}\begin{array}{r}\pe_{l_{\st 3}+1}\\
\vspace{5mm}\\ \\ \end{array}\!\!\!\Le\!\!\!\tilde\q_{l_{\st 3}+1}\er_{l_{\st
3}+2}\\ \!\!\!\pe_{l_{\st 3}+2}\q_{l_{\st 3}+2}\er_{l_{\st 3}+3}\\ \quad
\dt{2cm}\\ \hspace{15mm}\pe_{l_{\st 2}}\;\q_{l_{\st 2}}\!\!\!\R\!\!
\begin{array}{l}\\=\sx{2}{C}\!{}^{l_{\st 3}+1}_{l_{\st 2}}\\ \\0\end{array}\\
\hspace{34mm}\begin{array}{r}\pe_{l_{\st 2}+1}\\ \!\!\!
\sx{1}{C}\!{}^{l_{\st 2}+1}_{(l_1=m)}
\!\!=\\ \vspace{3mm}\end{array}\!\!\Le\!\!\!\tilde\q_{l_{\st 2}+1}\er_{l_{\st
2}+2}\\ \!\!\!\pe_{l_{\st 2}+2}\q_{l_{\st 2}+2}\er_{l_{\st 2}+3}\\ \ \ \dt{2cm}
\\ \hspace{17mm}\pe_m\ \q_m\!\!\!\R\vspace{1.3mm}\!\!\!\R\!\Rightarrow
\begin{array}{l}\!\!\!\!\!{}^{l_{n+1}\!+\!1=1}\st\hspace{4mm}l_n\hspace{12mm}l_2
\hspace{8mm}l_1=m\\
\!\!\underbrace{\Le\!\begin{picture}(45.00,53.00)
\multiput(45.00,15.50)(0,2){19}{\line(0,1){1}}
\multiput(30.00,33.50)(0,2){10}{\line(0,1){1}}
\put(15.00,53.50){\vector(0,1){3.00}}\put(30.00,53.50){\vector(0,1){3.00}}
\put(45.00,53.50){\vector(0,1){3.00}}\put(00.00,53.50){\vector(0,1){2.00}}
\put(-1.00,00.00){\vector(-1,0){4.00}}\put(-1.00,15.00){\vector(-1,0){4.00}}
\put(-1.00,18.00){\vector(-1,0){4.00}}\put(-1.00,33.00){\vector(-1,0){4.00}}
\put(-1.00,38.00){\vector(-1,0){4.00}}\put(-1.00,53.00){\vector(-1,0){4.00}}
\put(-6.00,02.00){\makebox(0,0)[cc]{$\st l_1=m$}}
\put(-6.00,13.00){\makebox(0,0)[cc]{$\st l_2+1$}}
\put(-6.00,20.00){\makebox(0,0)[cc]{$\st l_2$}}
\put(-6.00,31.00){\makebox(0,0)[cc]{$\st l_3+1$}}
\put(-6.00,36.00){\makebox(0,0)[cc]{$\st l_n$}}
\put(-6.00,51.00){\makebox(0,0)[cc]{$\st 1$}}
\put(0.00,0.00){\framebox(45.00,15.00)[cc]{$\sx{1}{B}$}}
\put(0.00,18.00){\framebox(30.00,15.00)[cc]{$\sx{2}{B}$}}
\multiput(0,35.00)(5,0){5}{$-\ $}
\put(0.00,38.00){\framebox(15.00,15.00)[cc]{$\sx{n}{B}$}}
\end{picture}\!\!\R}_{\sy{B}}\!,\!\!\end{array}$$
where
				$$
\tilde\q_{l_{k+1}+1}\!=\q_{l_{k+1}+1}\!-\!\pe_{l_{k+1}+1}\!\!
\sx{k+1}{B}\!\!\!{}_{l_{k+1}l_{k+1}}\!\er_{l_{k+1}+1}, k=1,2,..,n-1\leqno(2.7)
				$$
and $\sx{k+1}{B}\!\!\!{}_{l_{k+1}l_{k+1}}$ are the last diagonal elements of
submatrices $\sx{k+1}{B}$ which coincide with the last diagonal elements
of submatrices, inverse of well-posed submatrices $\sx{k+1}{C}^{l_{k+2}+1}_{k+1}$
separated by the method; and $n$ is the number of separated subspaces.

Elements $\sx{k}{B}^{}_{ij}$ are calculated [1] by the formulae:
$$\sx{k}{B}^{}_{ij}=\left\{\begin{array}{l}\omega^{}_i\PL^i_{\xi=j+1}\!\!\beta_
{\xi},\:\mbox{if}\:1\leq j<i,l_{k+1}+1\leq i\leq l_k,\\ 0\:\mbox{for all $\:
i\:$ from }\:j\!<\!i\!\leq\!l_k,\:\mbox{ if }\:\Lambda_j\!=\!0,\:\mbox{for any
$ \:j\:$ from }\:l_{k+1}+2\!\!\leq\!j\!\leq\!l_k,\\ 0\:\mbox{for all  $\:j\:$ from }\:1\!
\leq\!j\!<\!i,\:\mbox{ if }\sx{k}{\G}^{}_i=\!0,\:\mbox{ for any $\:i\:$ from }
\:l_{k+1}\!+\!1\!\leq\!i\!\leq\!l_k\!-\!1,\\
\omega^{}_i\PL^j_{\xi=i+1}\!\!\sx{k}{\hat\beta}_{\xi},
\mbox{ if } l_{k+1}+1\leq i<j\leq l_k,\\ 0\mbox{ for all  $i$  from }
l_{k+1}+1\leq i<j,\mbox{ if  }\sx{k}{\G}^{}_j=0,\\ 0\mbox{ for all  $j$  from }
 i<j\leq l_k,\mbox{ if  }\Lambda_i=0.\end{array}\right.\leqno(2.8)$$
Diagonal elements $\sx{k}{B}_{ii}$ of submatrices $\sx{k}{B}$ and
quantities $\omega^{}_i$ in (2.8) are calculated [1] by the formulae:
$$\left\{\begin{array}{l}\sx{k}{B}_{ii}=(\Lambda_{i+1}+\sx{k}{\G}_{i-1}-
\q_i)^{-1}\mbox{and }\omega^{}_i=\sx{k}{B}_{ii},\mbox{ if }\Lambda_i\ne
0\ne\sx{k}{\G}_i.\\ \sx{k}{B}_{ii}=0, \sx{k}{B}_{i-1i-1}=\sx{k}{\G}_{i-1}
\omega^{}_i,\sx{k}{B}_{i+1i+1}=\sx{k}{\G}\!{}^{-1}_i\mbox{and }\omega^{}_i
=(-\pe_i\er_i)^{-1},\mbox{ if }\Lambda_i=0.\\ \sx{k}{B}_{ii}=0,\sx{k}{B}_
{i-1i-1}=\Lambda^{-1}_i,\sx{k}{B}_{i+1i+1}=\Lambda_{i+1}\omega^{}_i\mbox{ and }
\omega^{}_i=(-\er_{i+1}\pe_{i+1})^{-1},\mbox{ if }
\sx{k}{\G}_i=0,\end{array}\right.\leqno(2.9)$$ 
sequences $\{\Lambda\}$ and $\{\sx{k}{\G}\}$ are computed by the formulae:
$$\left\{\!\!\!\begin{array}{l}\Lambda_{i+1}\!=\!\q_i\!-\!\pe_i\Lambda_i^{-1}\!\er_i,
\;\Lambda_2\!=\!\q_1,\;i\!=2,...,m,\mbox{ if }\Lambda_i\!\ne0\mbox{ for
all} \;2\!\leq\!i\!\leq\!m;\\ \mbox{if}\;\Lambda_i\!=\!0\;\mbox{for any
$i$ from}\;(\!2\!
\leq\!i\leq\!m\!),\mbox{then}\;\Lambda_{i+1}\!\!-\mbox{is undefined, but
}\;\Lambda_{i+2} \!\!=\!\q_{i+1};\end{array}\!\!\right.\leqno(2.10)$$ %
$$\left\{\!\!\!\begin{array}{l}\sx{k}{\G}_{i-1}\!=\!\q_i-\er_{i+1}\sx{k}{\G}\!{}_i^
{-1}\pe_{i+1},\sx{k}{\G}_{l_{\st k}-1}\!=\!\q_{l_{\st k}},i=l_k\!-1,l_k\!-2,
...,l_{k+1}\!+1,\;\mbox{if}\;\sx{k}{\G}_i\ne0;\\
\mbox{if}\;\sx{k}{\G}_i=\!  0\;\mbox{for any $i$
from}\;(l_{k+1}+1\!\leq\!i\leq\!l_k-1\!),\;\mbox{then}\;
\sx{k}{\G}_{i-1}\!\!-\mbox{is undefined,}\\ \mbox{but}\;\sx{k}{\G}_{i-2}\!=\!
\q_{i-1}.\end{array}\right.\leqno(2.11)$$
The structure elements $\beta_{\xi}$ and $\sx{k}{\hat\beta}_{\xi}$ which
determine the elements of submatrices $\sx{k}B$ and their products
$\prod\beta_{\xi}$ and $\prod \sx{k}{\hat\beta}_{\xi}$ are computed
[1] by the formulae:  %
$$\beta_i=\!\left\{\!\begin{array}{l}-\pe_i\Lambda^{-1}_i,\mbox{ if
}\Lambda_i \ne 0,\\-\pe_i,\mbox{ and
}\beta_{i+1}=-\pe_{i+1}\omega^{}_i,\\ \mbox{ if }
\Lambda_i=0;\!\end{array}\!\right.\sx{k}{\hat\beta}_{i+1}=\!\left\{\!
\begin{array}{l}-\er_{i+1}\sx{k}{\G}^{-1}_i,\mbox{ if }\sx{k}{\G}^{}_i\ne 0,
\\ -\er_{i+1},\mbox{ and }\sx{k}{\hat\beta}_i=-\er_i\omega^{}_i,
\\ \mbox{ if }\sx{k}{\G}^{}_i=0;\end{array}\right.\leqno(2.12)$$
$$\PL^i_{\xi=j+1}\!\!\!\beta_{\xi}\!=\!\left\{\!\!\!\begin{array}{l}\beta_i\!
\cdot\!\beta_{i-1}\!\cdots\beta_{j+1},\mbox{ if }j\!<\!i,\\1,\mbox{ if
}j\geq i;
\end{array}\right.\PL^j_{\xi=i+1}\!\!\!\sx{k}{\hat\beta}_{\xi}=\!\left\{
\!\!\!\begin{array}{l}\sx{k}{\hat\beta}_{i+1}\!\cdots\!\sx{k}{\hat\beta}_{j-1}
\!\cdot\!\sx{k}{\hat\beta}_j,\mbox{ if }i\!<\!j,\\1,\mbox{ if }i\geq j;
\end{array}\right.\leqno(2.13)$$

Proof. Let the system $\tilde C_3\bar X=\tilde Y$, according to the theorem
condition, be ill-posed but non-degenerate. Then its solution $X^+$ with the
properties given in the theorem does {\it theoretically} exist and it is
unique. Let us show that it can numerically be obtained by the method
$(2.1)\div(2.13)$ called in [1] the critical-component method. To this end,
we verify first that to the solution $X^+$ there corresponds the following
generalized $LDR$ [1] decomposition of the matrix $C_3$ (1.4):  %
$$C_3\!=\!LDR\!=\!\Le\!\!\!\Le\!1\\ \hspace{13mm}\ddots\\  \hspace{35mm}1\!\R
\\ \ddt{15}\\ \!\!\begin{array}{l}\st (\pe_{l_{\st 3}+1}\!\!\sx{3}{B}_{l_{\st
3}l_{\st 4}+1})\cdots(\pe_{l_{\st 3}+1}\!\!\sx{3}{B}_{l_{\st 3}l_{\st 3}})\\
\vspace{2.2mm} \\ \vspace{2mm} \end{array}\!\!\Le\!1\\
\hspace{13mm}\ddots\\ \hspace{33mm}1\!\R\\ \hspace{42mm}\!\!\begin{array}{l}\st
(\pe_{l_{\st 2}+1}\!\!\sx{2}{B}_{l_{\st 2}l_{\st 3}+1})\cdots(\pe_
{l_{\st 2}+1}\!\!\sx{2}{B}^{}_{l_{\st 2}l_{\st 2}})\\ \vspace{2mm} \\
\vspace{2mm} \end{array}\!\!
\Le\! 1\\ \hspace{5mm}\ddots\\ \hspace{12mm}1\!\R\!\!\R\!\!\times\leqno(2.14)$$
$$\times\!\!\!\Le\!\!\!\Le\!\!\!\q_1\ \er_2\\ \!\!\!\pe_2\ \q_2\ \er_3\\ \quad
\dt{1cm}\\ \hspace{7mm}\pe_{l_{\st n}}\q_{l_{\st n}}\!\!\!\R\!\!=
\sx{n}{C}\!{}^{(l_{\st n+1}+1=1)}_{l_{\st n}}\\ \hspace{13mm}\ddt{5}\\ \hspace{2cm}
\Le\!\!\!\tilde\q_{l_{\st 3}+1}\er_{l_{\st 3}+2}\\ \!\!\!\pe_{l_{\st 3}+2}\q_{l_{\st 3}+2}
\er_{l_{\st 3}+3}\\ \quad\cdots\cdots\cdots\\ \hspace{15mm}\pe_{l_{\st 2}}\;\q_
{l_{\st 2}}\!\!\!\R\!\!=\sx{2}{C}\!{}^{l_{\st 3}+1}_{l_{\st 2}}\\ \hspace{33mm}
\sx{1}{C}\!{}^{l_{\st 2}+1}_{(l_{\st 1}=m)}\!\!=\Le\!\!\!\tilde\q_{l_{\st 2}+1}\er_{l_{\st 2}+2}\\
\!\!\!\pe_{l_{\st 2}+2}\q_{l_{\st 2}+2}\er_{l_{\st 2}+3}\\ \ \ \dt{2cm}\\
\hspace{17mm}\pe_m\ \q_m\!\!\!\R\vspace{1.3mm}\!\!\!\R\!\!\!\Le\!\!\!\Le\!\!\! 1\\
\!\ddots\\ \hspace{5mm}1\!\!\!\R\!\!\!\begin{array}{l}\st (\sx{n}{B}_{1l_
{\st n}}\!\!\er_{l_{\st n}+1})\\ \qquad\vdots\\ \st(\sx{n}{B}_{l_{\st n}
l_{\st n}}\!\!\er_{l_{\st n}+1})\end{array}\\ \hspace{7mm}\ddt{5}\\
\hspace{9mm}\Le\!\!\! 1\\ \hspace{7mm}\ddots\\ \hspace{17mm} 1\!\!\R
\!\!\!\begin{array}{l}\st(\sx{2}{B}_{l_{\st 3}+1l_{\st 2}}\!\!\er_
{l_{\st 2}+1})\\ \qquad\vdots\\ \st(\sx{2}{B}^{}_{l_{\st 2}l_{\st 2}}\!\!
\er_{l_{\st 2}+1})\end{array}\\ \hspace{33mm} \Le\!\! 1\\ \hspace{8mm}\ddots\\
\hspace{2cm}1\!\!\!\R\vspace{1.5mm}\!\!\!\R\!,$$
where it is assumed that tridiagonal matrices $\sx{k}{C}\!{}^{l_{\st
k+1}+1}_{l_{\st k}}(k=1,2,...,n;l_1=m,l_{n+1}+1=1)$ are well-posed and their
first diagonal elements are denoted by
				$$
\tilde\q_{l_{k+1}+1}\!=\q_{l_{k+1}+1}\!-\!\pe_{l_{k+1}+1}\!\!
\sx{k+1}{B}\!\!\!{}_{l_{k+1}l_{k+1}}\!\er_{l_{k+1}+1},k=1,2,..,n-1,\leqno(2.14)'
				$$
where $\{\pe_i,\q_i,\er_i\}^{l_{\st k+1}}_{i=l_{\st k+2}+1}$ are elements of
the initial matrix $C_3$ (1.4), $\sx{k+1}{B}_{l_{\st k+1}j}(j=l_{\st
k+1},l_{\st k+1}-1,...,l_ {k+2}+1)$ are the last rows and
$\sx{k+1}{B}_{il_{\st k+1}}(i=l_{\st k+1},l_{\st k+1}-1,...,l_ {k+2}+1)$ are
the last columns of matrices, inverse of the matrices $\sx{k+1}{C}\!{}^{l_
{\st k+2}+1}_{l_{\st k+1}},$ computed in accordance with $\sx{k}{B}_{ij}$
(2.8) since they are elements of rectangular submatrices $\sx{k}{B}$
(2.8). From the assumptions for $C_3$ being nonsingular and for square
matrices $\sx{k}{C}\!{}^{l_{\st k+1}+1}_{l_ {\st k}}$ being well-posed it
follows that the $LDR$ decomposition (2.14) is unique and stable to errors
$(h,\varepsilon_1,\varepsilon_0)$.

And the matrix $B=C^{+}_3$ can uniquely be represented in the form
$(B=(E+\Omega)\sy{B})$:
$$B\!=\!(\sy{B}=\!\Le\!\!\!\Le\!\!\!\q_1\ \er_2\\ \!\!\!\pe_2\
\q_2\ \er_3\\ \quad\cdots\cdots\\ \hspace{7mm}\pe_{l_{\st n}}\q_{l_{\st n}}\!\!
\!\R^{-1}\!\!=(\sx{n}C\!{}^{(l_{\st n+1}+1=1)}_{l_{\st n}})^{-1}\\ \hspace{13mm}
\ddt{5}\\ \hspace{2cm}\Le\!\!\!\tilde\q_{l_{\st 3}+1}\er_{l_{\st 3}+2}\\
\!\!\!\pe_{l_{\st 3}+2}\q_{l_{\st 3}+2}\er_{l_{\st 3}+3}\\ \quad\cdots\cdots
\cdots\\ \hspace{15mm}\pe_{l_{\st 2}}\;\q_{l_{\st 2}}\!\!\!\R^{-1}\!\!=
(\sx{2}C\!{}^{l_{\st 3}+1}_{l_{\st 2}})^{-1}\\ \hspace{26mm}(\sx{1}{C}\!{}^{l_
{\st 2}+1}_{(l_{\st 1}=m)})^{-1}\!\!=\Le\!\!\!\tilde\q_{l_{\st 2}+1}\er_{l_{\st
2}+2}\\ \!\!\!\pe_{l_{\st 2}+2}\q_{l_{\st 2}+2}\er_{l_{\st 2}+3}\\ \ \ \cdots
\cdots\cdots\\ \hspace{17mm}\pe_m\ \q_m\!\!\!\R^{-1}\vspace{1.3mm}\!\!\!\R\!\!
\times\leqno(2.15)$$
$$\times\!\!\Le\!\!\!\Le\!1\\ \hspace{3cm}\ddots\\  \hspace{6cm}1\!\R\\
\ddt{15}\\ \!\!\begin{array}{l}\st (\!-\!\pe_{l_{\st 3}+1}\!\!\sx{3}{B}^{}_{l_
{\st 3}l_{\st 4}+1})\hspace{13mm}\cdots\hspace{12mm} \ (\!-\!\pe_{l_{\st
3}+1}\!\!\sx{3}{B}^{}_{l_{\st 3}l_{\st 3}})\\ \vspace{2.2mm} \\ \vspace{2mm}
\end{array}\Le\!1\\ \hspace{13mm}\ddots\\ \hspace{33mm}1\!\R\\
\!\!\begin{array}{l}\st (\!\pe_{l_{\st 2}+1}\!\!\sx{2}{B}^{}_
{l_{\st 2}l_{\st 3}+1}\!\pe_{l_{\st 3}+1}\!\!\sx{3}{B}^{}_{l_{\st 3}l_{\st
4}+1})\cdots(\!\pe_{l_{\st 2}+1}\!\!\sx{2}{B}^{}_{l_{\st 2}l_{\st 3}+1}\!
\pe_{l_{\st 3}+1}\!\!\sx{3}{B}^{}_{l_{\st 3}l_{\st 3}})(\!-\!\pe_{l_{\st 2}+1}
\!\!\sx{2}{B}^{}_{l_{\st 2}l_{\st 3}+1})\cdots(\!-\!\pe_{l_{\st 2}+1}\!\!\sx{2}
{B}^{}_{l_{\st 2}l_{\st 2}})\\ \vspace{2mm}\\ \vspace{2mm}\end{array}\!\!
\Le\! 1\\ \hspace{5mm}\ddots\\ \hspace{12mm}1\!\R\!\!\R\!)+$$
$$+(\Omega\!=\!\Le\!\!\!\Le\!\!\! 0\\ \!\ddots\\ \hspace{5mm}0\!\!\!\R\!
\!\!\begin{array}{l}\st (\!-\!\sx{n}{B}^{}_{1l_{\st n}}\!\!\er_{l_{\st n}+1})\;
{\displaystyle 0}\cdots{\displaystyle 0}
\;(\!\sx{n}{B}^{}_{1l_{\st n}}\!\!\er_{l_{\st n}+1}\!\sx{2}{B}^{}_{l_
{\st 3}+1l_{\st 2}}\!\!\er_{l_{\st 2}+1})\\ \qquad\vdots\hspace{1cm}\cdots
\hspace{1cm}\vdots\\ \st(\!-\!\sx{n}{B}^{}_{l_{\st n}l_{\st n}}\!\!
\er_{l_{\st n}+1})\;{\displaystyle 0}\cdots{\displaystyle 0}\;
\st(\!\sx{n}{B}^{}_{l_{\st n}l_{\st n}}\!\!\er_{l_{\st n}+1}\!\sx{2}{B}^{}_
{l_{\st 3}+1l_{\st 2}}\!\!\er_{l_{\st 2}+1})\end{array}\\
\hspace{7mm}\ddt{14}\\ \hspace{9mm}\Le\!\!\! 0\\ \hspace{8mm}\ddots\\
\hspace{2cm}0\!\!\R
\!\!\!\begin{array}{l}\st(\!-\!\sx{2}{B}^{}_{l_{\st 3}+1l_{\st 2}}\!\!\er_
{l_{\st 2}+1})\\ \qquad\vdots\\ \st(\!-\!\sx{2}{B}^{}_{l_{\st 2}l_{\st 2}}\!\!
\er_{l_{\st 2}+1})\end{array}\\ \hspace{38mm} \Le\!\!0\\ \hspace{12mm}\ddots\\
\hspace{33mm}0\!\!\!\R\vspace{1.5mm}\!\!\!\R\!)\cdot\!\sy{B}.$$
Schematically, the matrix $\sy{B}$ can be represented as follows:
				$$
\begin{array}{l}\!\!\!\!\!{}^{l_{n+1}\!+\!1=1}\st\hspace{4mm}l_n
\hspace{12mm}l_2\hspace{8mm}l_1=m\\
\Le\!\begin{picture}(45.00,53.00)
\multiput(45.00,15.50)(0,2){19}{\line(0,1){1}}
\multiput(30.00,33.50)(0,2){10}{\line(0,1){1}}
\put(15.00,53.50){\vector(0,1){3.00}}\put(30.00,53.50){\vector(0,1){3.00}}
\put(45.00,53.50){\vector(0,1){3.00}}\put(00.00,53.50){\vector(0,1){2.00}}
\put(-1.00,00.00){\vector(-1,0){4.00}}\put(-1.00,15.00){\vector(-1,0){4.00}}
\put(-1.00,18.00){\vector(-1,0){4.00}}\put(-1.00,33.00){\vector(-1,0){4.00}}
\put(-1.00,38.00){\vector(-1,0){4.00}}\put(-1.00,53.00){\vector(-1,0){4.00}}
\put(-6.00,02.00){\makebox(0,0)[cc]{$\st l_1=m$}}
\put(-6.00,13.00){\makebox(0,0)[cc]{$\st l_2+1$}}
\put(-6.00,20.00){\makebox(0,0)[cc]{$\st l_2$}}
\put(-6.00,31.00){\makebox(0,0)[cc]{$\st l_3+1$}}
\put(-6.00,36.00){\makebox(0,0)[cc]{$\st l_n$}}
\put(-6.00,51.00){\makebox(0,0)[cc]{$\st 1$}}
\put(0.00,0.00){\framebox(45.00,15.00)[cc]{$\sx{1}{B}$}}
\put(0.00,18.00){\framebox(30.00,15.00)[cc]{$\sx{2}{B}$}}
\multiput(0,35.00)(5,0){5}{$-\ $}
\put(0.00,38.00){\framebox(15.00,15.00)[cc]{$\sx{n}{B}$}}
\end{picture}\!\!\R=\sy{B}.\!\!\end{array}\leqno(2.15)'$$
Representation (2.15) is easily established by a direct verification of the
matrix equalities $C_3B=E=BC_3$, with representations $(2.8)\div(2.13)$
taken into account for elements $\sx{k}B_{ij}$ of the matrix $\sy{B}$
and decompositions of $C_3$ and $B$ given by (2.14) and (2.15).

Now using representation (2.15) for $B$, we obtain components of the vector
$X^+$ in form (2.2) and (2.5). From (2.15) it follows that $X^+$ can be
written in the form %
$$X^+=(E+\Omega)\sy{X}=\sy{X}+\Omega\sy{X},\leqno(2.16)$$
where the vector $\sy{X}$ looks as follows
$$\sy{X}=(\sy{B}\equiv\sy{C}_3\!{}^{-1})Y\leqno(2.17)$$
and is a unique, stable to errors $(h,\delta)$ and $(\varepsilon_
1,\varepsilon_0)$, solution of the well-posed system of linear algebraic
equations
$$(\sy{C}_3=\!\!\Le\!\!\!\Le\!\!\!\q_1\ \er_2\\
\!\!\!\pe_2\ \q_2\ \er_3\\ \quad\dt{1cm}\\ \hspace{7mm}\pe_{l_{\st
n}}\q_{l_{\st n}}\!\!\!\R\!\!  \begin{array}{l}\\=\sx{n}{C}\!{}^{(l_{\st
n+1}+1=1)}_{l_{\st n}}\\ \\0\end{array}\\ \hspace{13mm}\ddt{5}\\
\hspace{1cm}\begin{array}{r}\pe_{l_{\st 3}+1}\\ \vspace{5mm}\\ \\
\end{array}\!\!\!\Le\!\!\!\tilde\q_{l_{\st 3}+1}\er_{l_{\st 3}+2}\\
\!\!\!\pe_{l_{\st 3}+2}\q_{l_{\st 3}+2}\er_{l_{\st 3}+3}\\ \quad \dt{2cm}\\
\hspace{15mm}\pe_{l_{\st 2}}\;\q_{l_{\st 2}}\!\!\!\R\!\!
\begin{array}{l}\\=\sx{2}{C}\!{}^{l_{\st 3}+1}_{l_{\st 2}}\\ \\0\end{array}\\
\hspace{34mm}\begin{array}{r}\pe_{l_{\st 2}+1}\\ \!\!\!
\sx{1}{C}\!{}^{l_{\st 2}+1}_{(l_{\st 1}=m)}
\!\!=\\ \vspace{3mm}\end{array}\!\!\Le\!\!\!\tilde\q_{l_{\st 2}+1}\er_{l_{\st
2}+2}\\ \!\!\!\pe_{l_{\st 2}+2}\q_{l_{\st 2}+2}\er_{l_{\st 2}+3}\\ \ \ \dt{2cm}
\\ \hspace{17mm}\pe_m\ \q_m\!\!\!\R\vspace{1.3mm}\!\!\!\R\!)\!\cdot\!\!
\underbrace{\Le\!\!\!\vspace{1.5mm}
\Le\sx{n}{x}^{}_1\\ \vdots\\ \sx{n}{x}^{}_{l_{\st n}}\!\!\!\R\!\!\!\\
\ \vdots\\ \!\!\!
\Le\!\!\sx{2}{x}^{}_{l_{\st 3}+1}\!\!\!\\ \vdots\\ \sx{2}{x}^{}_{l_{\st 2}}\!\!\!\R
\!\!\!\\ \!\!\!
\Le\!\!\sx{1}{x}^{}_{l_{\st 2}+1}\!\!\!\\ \vdots\\ \sx{1}{x}^{}_m\!\!\!\R
\vspace{1.3mm}\!\!\!\R}_{\sy{X}}\!\!=\!\!\underbrace{\Le\!\!\!\vspace{1.8mm}
\Le {y}^{}_1\\ \vdots\\ {y}^{}_{l_{\st n}}\!\!\!\R\!\!\!\\
\ \vdots\\ \!\!\!
\Le\!\! {y}^{}_{l_{\st 3}+1}\!\!\!\\ \vdots\\ {y}^{}_{l_{\st 2}}\!\!\!\R
\!\!\!\\ \!\!\!
\Le\!\! {y}^{}_{l_{\st 2}+1}\!\!\!\\ \vdots\\ {y}^{}_m\!\!\!\R
\vspace{1.3mm}\!\!\!\R}_{Y}\!,\leqno(2.18)$$
which differs from the initial system $C_3X=Y$ by the change of the
corresponding off-diagonal elements to zeros and of diagonal elements $q$ to
elements $\tilde q$ calculated by formulae $(2.14)'$. Here the vector
$\sy{X}$ includes components given by sums (2.2), which follows from
the representation $\sy{B}$ (2.15) and (2.8).

For the matrix $\Omega$ (2.15) we can write the following decomposition
$$\Omega\!=\!\Le\!\!\!\Le\!\!\! 0\\
\!\ddots\\ \hspace{5mm}0\!\!\!\R\!\!\!\begin{array}{l}\st (\!-\!
\sx{n}{B}^{}_{1l_{\st n}}\!\!\er_{l_{\st n}+1})\\ \qquad\vdots\\ \st(\!-\!\sx{n}
{B}^{}_{l_{\st n}l_{\st n}}\!\!\er_{l_{\st n}+1})\end{array}\\ \hspace{7mm}
\ddt{5}\\ \hspace{1cm}\Le\!\!\!0\\ \hspace{7mm}\ddots\\ \hspace{17mm}0\!\!\R
\!\!\!\begin{array}{l}\st(\!-\!\sx{2}{B}^{}_{l_{\st 3}+1l_{\st 2}}\!\!\er_
{l_{\st 2}+1})\\ \qquad\vdots\\ \st(\!-\!\sx{2}{B}^{}_{l_{\st 2}l_{\st 2}}\!\!
\er_{l_{\st 2}+1})\end{array}\\ \hspace{34mm} \Le\!\!0\\ \hspace{8mm}\ddots\\
\hspace{2cm}0\!\!\!\R\vspace{1.5mm}\!\!\!\R
\!\!\cdots\!\!
\Le\!\!\!\Le\!\!\!1\\ \!\ddots\\ \hspace{5mm}1\!\!\!\R\\ \hspace{7mm}
\ddt{3}\\ \hspace{1cm}\Le\!\!\! 1\\ \!\ddots\\ \hspace{5mm}1\!\!\!\R
\!\!\!\begin{array}{l}\st(\!-\!\sx{2}{B}^{}_{l_{\st 3}+1l_{\st 2}}\!\!\er_
{l_{\st 2}+1})\\ \qquad\vdots\\ \st(\!-\!\sx{2}{B}^{}_{l_{\st 2}l_{\st 2}}\!\!
\er_{l_{\st 2}+1})\!\!\end{array}\!\!\\ \hspace{22mm}\Le\!\!1\\ \hspace{6mm}
\ddots\\ \hspace{17mm}1\!\!\!\R\vspace{1.5mm}\!\!\!\R.\leqno(2.19)$$
Then for components of the vector  $\Omega\cdot\!\!\sy{X}$, denoted as
vector $\phi$, we obtain the explicit form $$
\Omega\cdot\!\!\sy{X}=[(-\sx{n}B\!{}_{1l_{\st n}}\er_{l_{\st n}+1})x^+_
{l_{\st n}+1},...,(-\sx{n}B\!{}_{l_{\st n}l_{\st n}}\er_{l_{\st n}+1})x^+_
{l_{\st n}+1};...;(-\sx{2}B\!{}_{l_{\st 3}+1l_{\st 2}}\er_{l_{\st 2}+1})x^+_
{l_{\st 2}+1},...\leqno(2.20)$$
{}\hfill $...,(-\sx{2}B\!{}_{l_{\st 2}l_{\st 2}}\er_{l_
{\st 2}+1})x^+_{l_{\st 2}+1};0,...,0]^T=\phi.$\\
Consequently, components of the vector $\phi$ are also calculated by
formulae (2.2).

As a result, we have established that if $X^+$ is a normal (pseudo)solution
of the system $C_3X=Y$ with the properties given by the theorem, then to it
there correspond consistent with it decompositions (2.14) and (2.15) for
matrix $C_3$ and its (pseudo)inverse matrix $B\equiv C^+_3$. In this case,
representations for $X^+=\sy{X}+\Omega\cdot\!\!\sy{X}$ (2.16) and
$C^+_3=\sy{B}+ \Omega\cdot\!\!\sy{B}$ (2.15) being consistent with each
other (since they are calculated with the same matrix $\Omega$ (2.15)) are
unique and stable to small errors $h,\delta)$ and
$(\varepsilon_1,\varepsilon_0)$ in view of decompositions of $C_3$ (2.14)
and $C^+_3$ (2.15) being unique. Stability is a consequence of the matrix
$\sy{C}_3$ (2.18) being well-posed.

Now let us show that if the numerical solution of the system $\tilde C_3\bar
X=\tilde Y$ is obtained by the method\footnote[0]{As we can see, this method
includes the algorithm and criterion $(2.3)\div(2.4)$ of separation of
well-posed subspaces and, respectively, the procedure of numerical finding
of $X^+$. It is to be kept in mind that the quantities $\Phi^{}_j$
(2.4) obey the inequalities $|\Phi^{}_j|\leq|\Delta^{}_j|,$ where
$\Delta^{}_j $ is a discrepancy. As a matter of fact,
$|\Phi^{}_j|=||\y_j|-|\y_j-[\y_j-(\pe_j\sx{k}x\!{}_{j-1}+\q_j\sx{k}x\!{}_j+\
r_{j+1}\sx{k}x\!{}_{j+1})]||=||\y_j|-|\y_j-\Delta^{}_j||\leq|\Delta^{}_j|,
$since $||\y_j|-|\y_j-\Delta^{}_j||\leq(|\y_j-(\y_j+\Delta^{}_j)|=|\Delta^{}_j|).$}
$(2.1)\div(2.13),$ then it is minimal in norm and provides a
minimum of the discrepancy norm. Indeed, let $X^+$ is determined in the form
$(2.1)\div(2.13)$. Then from (2.5) it follows that the vector $X^+$ can be
represented as a sum of two vectors,

				$$
			X^+=\sy{X}+\phi,\leqno(2.21)
				$$
whose components are determined according to (2.2). The vectors
$\sy{X}$ and $\phi$ consist of $n$ subvectors of proper
dimensions, i.e.,
$$\begin{array}{l}\!\sy{X}\!=\![(\sx{n}{x}_1\!,...,\sx{n}{x}_{l_{\st n}}),...,\!
(\sx{k}{x}_{l_{\st k+1}+1}\!,...,\sx{k}{x}_{l_{\st k}}),...,\!(\sx{2}{x}_{l_{\st
3}+1},...,\sx{2}{x}_{l_{\st 2}}),...,\!(\sx{1}{x}_{l_{\st 2}+1}=x^+_{l_{\st 2}
+1},...,\sx{1}{x}_m\!=\!x^+_m\!)]^T\!,\\ \!\phi\!=\![(\sx{n}{\phi}_1,...,\sx{n}{\phi}_{l_
{\st n}}),...,(\sx{k}{\phi}_{l_{\st k+1}+1},...,\sx{k}{\phi}_{l_{\st k}}),...,
(\sx{2}{\phi}_{l_{\st 3}+1},...,\sx{2}{\phi}_{l_{\st 2}}),...,(\sx{1}{\phi}_
{l_{\st 2}+1}=0,...,\sx{1}{\phi}_m=0)]^T,\end{array}$$
corresponding to well-posed $n$-subspaces that are separated in accordance
with the criterion (2.3), (2.4).

>From $(2.1)\div(2.5)$ we have that to the solution $X^+$ there corresponds
the decomposition of $B=C^{+}$ of the form (2.15), which results in the
representation for $C_3$ (1.4) of form (2.14). Then, as mentioned above,
the solution $X^+$ is written in form (2.16), where $\sy{X}$ is in a unique
way represented in form (2.17), (2.18). Owing to system (2.18) being
well-posed, which results from criterion (2.3) and (2.4), the vector
$\sy{X}$ is unique, obeys the condition $\min||\sy{X}||$, and is stable
to small errors $(h,\delta)$ and
$(\varepsilon_1,\varepsilon_0)$. From the uniqueness of matrix $\Omega$
(2.15) that contains the last columns of matrices, inverse of the well-posed
matrices $\sx{k}{C}\!{}^{l_{k+1}+1}_ {l_k}$, and from (2.16) and
(2.15) it follows that $X^+$ and $(B\equiv C^+_3)$ are unique and minimal
in norm.

Let us now show that the vector $X^+$ determined by formulae
$(2.1)\div(2.5)$ satisfies the condition of minimum of the discrepancy norm
$(\min||\tilde C_3X^+-\tilde Y||)$. Taking advantage of the representation of
$X^+$ (2.15), we get $$ ||\tilde C_3X^+-\tilde Y||=||\tilde
C_3(E+\Omega)\sy{X}-\tilde Y||=||\sy{C}_3 \sy{X}-\tilde Y||, $$ where
$\sy{X}=(\sy{B}=\sy{C}\!{}^{-1}_3)\tilde Y,\sy{C}\!{}^{-1}_3$, and $\sy{B}$
are, respectively, defined by (2.18) and (2.15). Owing to the system
(2.16) being well-posed, the minimum $\min||\sy{C}_3\sy{X}-\tilde Y||$ and,
consequently, $\min||\tilde C_3X^+-\tilde Y||$ are attainable. So, the
theorem is proved.

{\bf Corollary.} The norm of discrepancy $||\tilde C_3X^+-\tilde Y||$
complies with the following estimate:  $$ \left\{\!\!\begin{array}{l}||\tilde
C_3X^+\!-\tilde Y||_{\infty}\!\leq\varepsilon_1
\tau\rho\gamma\max\limits_{1\leq i\leq m}|\tilde\y_i|+\Delta,\mbox{ where }
\tau=\max\limits_{2\leq i\leq m}(|\tilde\pe_i|,|\tilde\er_i|),
\rho=\max\limits_{i,j;k}(|\sx{k}{B}_{ij}|),\\ \gamma=\sqrt{\sum\limits^n_{k=1}
l_k(l_k-l_{k+1})},l_1=m,l_{n+1}=0,0\leq\Delta\leq h||X^+||+\delta.
\end{array}\right.\leqno(2.22)
				$$

Proof. In view of all said above, we have
				$$
||\tilde C_3X^+-\tilde Y||=||\tilde C_3(E+\Omega)\sy{X}-\tilde Y||=||\sy{C}_3
\sy{X}-\tilde Y||.
				$$
Since the system $\sy{C}_3\sy{X}=\tilde Y$ is well-posed, the Euclidean norm
of errors $||\sy{C}_3\sy{X}-\tilde Y||_E$ can be estimated by using the
known results [9]:
$$
||\tilde C_3X^+-\tilde Y||_E=||\sy{C}_3\sy{X}-\tilde Y||_E\leq4f(m)
\varepsilon_1||\sy{C}_3||_E||\sy{X}||_E.\leqno(2.23)
$$
However, in the case of the method considered above, this estimate turns out
to be excessive. Actually, performing obvious transformations and making
use of the definition of matrix norms consistent with the corresponding
vector norms,\linebreak\newpage
\noindent we get
$$
||\sy{C}_3\sy{X}-\tilde Y||=||\sy{C}_3\sy{B}\tilde Y-\tilde Y||=||(\sy{C}_3
\sy{B}-E)\tilde Y||\leq||\sy{C}_3\sy{B}-E||\cdot||\tilde Y||.\leqno(2.24)
$$
Next, we estimate the norm of matrix discrepancy $||\sy{C}_3\sy{B}-E||$,
using the explicit form of $\sy{C}_3$ (2.18) and $\sy{B}$ (2.15), as well as
the condition for the matrix $\sy{C}_3$ being well-posed. Taking account of
the explicit form of elements of the matrix $(\sy{C}_3\sy{B}-E)$ and
introducing the notation $\nu_{ij}(\nu=(\nu_{ij})=\sy{C}_3\sy{B}-E)$ for
them, we can write the system of scalar identities $$
\left\{\begin{array}{l}\tilde\pe_i\sx{k}{B}_{i-1j}+\tilde\q_i\sx{k}{B}_{ij}+
\tilde\er_{i+1}\sx{k}{B}_{i+1j}\equiv\nu_{ij},l_{k+1}+1\leq i<j\leq l_k,\\
\tilde\pe_i\sx{k}{B}_{i-1i}+\tilde\q_i\sx{k}{B}_{ii}+\tilde\er_{i+1}\sx{k}{B}_
{i+1i}\equiv1+\nu_{ii},l_{k+1}+1\leq(i=j)\leq l_k,\\ \tilde\pe_i\sx{k}{B}_{i-1j}+
\tilde\q_i\sx{k}{B}_{ij}+\tilde\er_{i+1}\sx{k}{B}_{i+1j}\equiv\nu_{ij},1\leq
j<i\leq l_k.\end{array}\right.\leqno(2.25)
$$
Hereafter, $k=1,2,...,n;l_1=m,l_{n+1}=0.$ Utilizing the representations for
$\sx{k}{B}_{ij}\ (2.8)\div(2.13)$, we write the system of identities (2.25)
either in the form
$$
\left\{\begin{array}{l}[\Lambda_{i+1}(-\Lambda^{-1}_{i+1}\tilde\er_{i+1})+
\tilde\er_{i+1}]\sx{k}{B}_{i+1i+1}\PL^j_{\xi=i+2}\sx{k}{\hat\beta}_{\xi}\equiv\nu_
{ij},l_{k+1}+1\leq i<j\leq l_k,\\ (\Lambda_{i+1}+\sx{k}{\G}_{i-1}-\tilde\q_i)
\sx{k}{B}_{ii}\equiv1+\nu_{ii},l_{k+1}+1\leq(i=j)\leq l_k,\\ \bigl[\sx{k}{\G}_{i-1}
(-\sx{k}{\G}\!{}^{-1}_{i-1}\tilde\pe_i)+\tilde\pe_i\bigr]\sx{k}{B}_{i-1i-1}
\PL^i_{\xi=j+1}\sx{k}{\beta}_{\xi}\equiv\nu_{ij},1\leq j<i\leq l_k,\end{array}
\right.\leqno(2.26)
$$
or in the form
$$
\left\{\begin{array}{l}(1-\Lambda_{i+1}\Lambda^{-1}_{i+1})\tilde\er_{i+1}
\sx{k}{B}_{i+1j}\equiv\nu_{ij},l_{k+1}+1\leq i<j\leq l_k,\\ (\Lambda_{i+1}+
\sx{k}{\G}_{i-1}-\tilde\q_i)\sx{k}{B}_{ii}\equiv1+\nu_{ii},l_{k+1}+1\leq(i=j)
\leq l_k,\\ (1-\sx{k}{\G}_{i-1}\sx{k}{\G}\!{}^{-1}_{i-1})\tilde\pe_i\sx{k}{B}_
{i-1j}\equiv\nu_{ij},1\leq j<i\leq l_k.\end{array}\right.\leqno(2.27)
$$
Let us now estimate $(2.27)_{1)},(2.27)_{2)}$ and $(2.27)_{3)}$; we have
$$
\left\{\!\!\begin{array}{l}\!|(1-\Lambda_{i+1}\Lambda^{-1}_{i+1})\tilde\er_
{i+1}\!\sx{k}{B}_{i+1j}\equiv\!\nu_{ij}|\leq|1-(1\pm\varepsilon_1)|\max\limits_
{1\leq i\leq\ m-1}|\tilde\er_{i+1}|\max\limits_{1\leq i<j\leq m-1}|\!\sx{k}{B}_
{i+1j}\!|,\\ |(\Lambda_{i+1}+\sx{k}{\G}_{i-1}-\tilde\q_i)\!\sx{k}{B}_{ii}-1
\equiv\nu_{ii}|\leq|(1\pm\varepsilon_1)-1|\max\limits_{1\leq i\leq m}|\sx{k}{B}_
{ii}|,\\ |(1-\sx{k}{\G}_{i-1}\sx{k}{\G}\!{}^{-1}_{i-1})\tilde\pe_i\!\sx{k}{B}_
{i-1j}\equiv\nu_{ij}|\leq|1-(1\pm\varepsilon_1)|\max\limits_{2\leq i\leq l_{\st
k}}|\tilde\pe_i|\max\limits_{2\leq j<i\leq l_{\st k}}|\sx{k}{B}_{i-1j}|.\!\!\!
\end{array}\!\!\right.\leqno(2.28)
$$
With estimates (2.28), we obtain
$||\tilde C_3X^+-\tilde Y||_{\infty}=||(\sy{C}_3\sy{B}-E)\tilde Y||_{\infty}
\leq||\sy{C}_3\sy{B}-E||_M||\tilde Y||_{\infty}\leq\sqrt{\sum\limits_{k=1}^nl_k
(l_k-l_{k+1})}\max\limits_{i,j}|\nu_{ij}|||\tilde Y||_{\infty}\leq\varepsilon_1
\tau\rho\gamma\max\limits_{1\leq i\leq m}|\tilde\y_i|,$
where $\tau,\rho,\gamma$ are defined by (2.22). Here we took advantage of the
condition of consistency of vector norms $||\tilde C_3X^+-\tilde
Y||_{\infty}=\max \limits_{1\leq i\leq m}|(\tilde C_3X^+-\tilde Y)_i|$ and
$||\tilde Y||_{\infty}= \max\limits_{1\leq i\leq m}|\tilde\y_i|$ with the
$M$-norm of the matrix $(\sy{C}_3\sy{B} -E)$, i.e.
$||\sy{C}_3\sy{B}-E||_M=\sqrt{m^2}\max\limits_{1\leq i,j\leq m}
|(\sy{C}_3\sy{B}-E)_{ij}|$. The validity of inequality (2.22) is established.
Since the Euclidean norm of the matrix $||\sy{C}_3\sy{B}-E||_E$ is consistent
only with the vector Euclidean norm $||\tilde Y||_E$, instead of (2.22), one
can, by analogous arguments, obtain the estimate $||\tilde C_3X^+-\tilde Y||_E=||(\sy{C}_3\sy{B}-E)\tilde Y||_E\leq
\varepsilon_1\tau\rho\gamma||\tilde Y||_E,$ where $\tau,\rho,\gamma$
are defined by (2.22).

{\bf Remark 1.} To save the volume of publication, we do not present the
method of solution of system (1.3) with the two-diagonal matrix $C_2$ (1.4).
It is expounded in detail in ref.[3] and it is shown there that it results
from the method $(2.1)\div(2.13)$.

The estimate (2.22) for system (1.3) acquires the following
form:
$$
\left\{\begin{array}{l}||\tilde C_2\hat X^+-\tilde{\hat Y}||\leq ||\sy{C}_2\sy{B}-
E||\cdot||\tilde Y||\leq\varepsilon_1\hat\tau\hat\rho\hat\gamma\max\limits_{1
\leq i\leq m}|\tilde\y_i|+\Delta,\mbox{where
}\hat\tau=\max\limits_{2\leq\i\leq m}(| \tilde\er_i|),\\
\hat\rho=\max\limits_{i,j;k}(|\sx{k}{B}_{ij}|),\hat\gamma=
\sqrt{1/2\sum\limits^n_{k=1}(l_k-l_{k+1})},l_1=m,l_{n+1}=0,0\leq\Delta\leq
h||\hat X^+||+\delta.\end{array}\right.
\leqno(2.29)
$$
Here $\sx{k}{B}_{ij}$ are elements of upper triangular matrices, inverse
of well-posed two-diagonal matrices $[\sx{k}{C}_2]^{l_{\st
k+1}+1}_{l_{\st k}}$.

{\bf Remark 2.} Note that due to orthogonality of matrices $P$ and $Q$ in
transformations (1.5), the following estimates take place for the norms of
discrepancy $||\tilde AZ^+-\tilde F||$: $$
\left\{\begin{array}{l}||\tilde AZ^+-\tilde F||\leq\varepsilon_1\tilde\tau
\tilde\rho\tilde\gamma\max\limits_{1\leq i\leq m}|\tilde\y_i|+\Delta,\mbox{
if }A=A^T,\\ ||\tilde AZ^+-\tilde
F||\leq\varepsilon_1\hat{\tilde\tau}\hat{\tilde\rho}
\hat{\tilde\gamma}\max\limits_{1\leq i\leq m}|\tilde\y_i|+\Delta,\mbox{ if
}A\ne A^T,\end{array}\right.\leqno(2.30) $$ where
$\tilde\tau,\tilde\rho,\tilde\gamma$ and $\hat{\tilde\tau},\hat{\tilde\rho},
\hat{\tilde\gamma}$ are defined in analogy with (2.22) and (2.29),
$0\leq\Delta\leq h|| Z^+||+\delta$.

{\bf Remark 3.} The above estimates (2.22), (2.29) and (2.30) can also be
used for problems of inversion i.e., $C_3C^+_3=E,C_2C^+_2=E,AA^+=E $:
the matrices $C^+_3,C^+_2$ and $A^+$ are to be obtained by solving the matrix
system of equations
$$C_3C^+_3=E, C_2C^+_2=E\mbox{ and }AA^+=E.$$
by the critical-component method.
In the case when systems (1.1), (1.2) and (1.3) are ill-posed, one should not
take, as $C^+_3,C^+_2$ and $A^+$, the corresponding matrices obtained by the
critical-component method in solving these systems of equations with a
given right-hand side. The reason is that the norms of matrices
$C^+_3,C^+_2$ и $A^+$ are consistent with the norms of concrete vectors
$X^+_3,\hat X^+_2,Z^+,\tilde Y,\tilde{\hat Y}$ and $\tilde F$.

{\bf Remark 4.} The theorem is formulated under the assumption $\det\tilde
C\ne0$. Let us remove this restriction. The critical-component method does
not explicitly use the quantity $\det\tilde C.$ Rather, it is based on the
processes (2.10) and (2.11) for computing elements of $m$-dimensional vectors
$\{\Lambda,\G\}$. As established in ref. [14], if $\det\tilde C=0$, then
components of these vectors get into one of the following three situations:
either $\Lambda_{m+1}=0$ and $\G_0=0$, or $\Lambda_i=0$ and $\G_i=0$, or
[($\Lambda_i=0$ and $\Lambda_i\Lambda_{i+1}=0$) or ($\G_i=0$ and
$\G_i\G_{i-1}=0$)]. In this case we replace some zero quantities by the
quantity $o(\varepsilon_1)$. This does not essentially impair the quality
of solution, since such perturbations can already be present in these
quantities. Consequently, one may consider the critical-component method
to be applicable for any value of $\det\tilde C,$ including
$\det\tilde C=0$.

\subsubsection{Results of numerical; experiments and their analysis}

\abzats
In this section, we discuss the results of numerical experiments performed
in the computer arithmetic with double accuracy
($\varepsilon_1=2^{-52}\approx 2,2\cdot10^{-16}$) for computing basic
numerical characteristics of the solutions $X$ of systems $WX=Y,\ W:\ C_2;
C_3; A\ne A^T; A=A^T.$ Let us first explain the notation and abbreviations
adopted in Tables $1\div12:\ \delta^{(m)}_M$ --- the relative
error of $\tilde X^ {(m)}$ --- the obtained numerical solution of system
$W^{(m)}X^{(m)}=Y^{(m)}\ (W^{(m)}:\ C^{(m)}_2; C^{(m)}_3; A^{(m)}\ne
(A^{(m)})^T; A^{(m)}=(A^{(m)})^T$ --- the above indicated types of matrices,
$X^{(m)}$ --- the exact solution, $m$ --- the order of the system under
consideration);
$\mu(W^{(m)})={\rm cond}(W^{(m)})$ --- the condition number of $W^{(m)};\
t^{(m)} (sec.)={\rm com. time (sec.)}$ --- the time of computing solutions
$\tilde X^{(m)},\ \delta^{(m)}_L,\ \delta^{(m)}_R$ --- the lower and upper
bounds $\delta^{(m)}_M$ i.e.\footnote[0]{The left-hand side of inequality
(3.1) is a property [10] of the norm $||\cdot||$, and the right-hand side is
obtained by using the exact solution $X=W^{-1}Y$ and equality
$W^{-1}W =E$. We have $(\delta_M=||\tilde X-X||/||X||)=||\tilde
X-W^{-1}Y||/||X||= ||W^{-1}(W\tilde X-Y)||/||X||\leq
(||W^{-1}||\cdot||W\tilde X-Y||/||X||= \delta_R).$ Note that in practice,
inequalities (3.1) can be broken (see, for instance, Table 7). This occurs
when calculating $W\tilde X-Y$. In this case, the solution
 $\tilde X$ can be surely considered acceptable.} $$
(\delta^{(m)}_L\!\!=\!\frac{\st\bigl| ||\tilde X^{(m)}\!||-||X^{(m)}\!||
\bigr|}{\st||X^{(m)}||})\!\leq\!(\delta^{(m)}_M\!\!=\!\frac{\st||\tilde
X^{(m)}\!-X^{(m)}\!||_E}{\st||X^{(m)}||})\!\leq\!(\frac{\st||(W^{(m)})^{-1}\!
||\cdot||W^{(m)}\tilde X^{(m)}\!-Y^{(m)}\!||_E}{\st||X^{(m)}||}\!=\!\delta^{(m)}_
R\!),\leqno(3.1)
$$
where $||\tilde X^{(m)}||_E,||X^{(m)}||_E$ are norms of approximate and exact
solutions; $\delta_L\!=\!\frac1N_j\sum\limits^{N_j}_{l=1}(\delta^{(m)}_L)_l,\\
\delta_M\!=\frac1N_j\sum\limits^{N_j}_{l=1}(\delta^{(m)}_M)_l,
\delta_R=\frac1N_j\sum\limits^{N_j}_{l=1}(\delta^{(m)}_R)_l,
\delta_{\tilde X}^*=\frac1N_j\sum\limits^{N_j}_{l=1}(||\tilde X^{(m)}||)_l,
\delta_X=\frac1N_j\sum\limits^{N_j}_{l=1}(||X^{(m)}||)_l,\\
\bar\mu(W)=\frac1N_j\sum\limits^{N_j}_{l=1}(\mu(W^{(m)}))_l,
t(sec.)=\frac1N_j\sum\limits^{N_j}_{l=1}(t^{(m)})_l$ are arithmetic means of
the characteristics listed above, $\bar
N_i=\sum\limits^s_{j=1}N_j$, where $N_j$ is the number of examples of a
given type, $s$ -- the number of examples in  a Table,
$i$ -- the number of a Table; {\tt MCS} and {\tt MCC} --- our programs {\tt
DCSOL} (access through  www {\tt http://cv.jinr.ru/lcta/sap/ lib/f499.f})
from library {\tt LIBJINR} [17] (algorithms of the critical-component method
$[1\div3]$); {\tt GS} --- programs {\tt DBEQN} and {\tt DEQN} from library
{\tt CERNLIB} [21] (a modified algorithm of the Gauss exclusion method);
{\tt OSM} --- program {\tt DTSYS} from library {\tt LIBJINR} [19]
(algorithm of nonmonotone orthogonal run); {\tt QR} --- programs {\tt
F01AXF} from library {\tt NAGLIB} [20] (algorithms of {\tt QR} --- method);
{\tt SVD} --- subprogram-function {\tt PSOL} from library {\tt LINA} [7]
(algorithm of the singular-expansion method with the use of exhaustion); {\tt
TRM} --- subprogram {\tt SLAY} from library {\tt LIBJINR} [19]  (algorithms
of the Tikhonov regularization method). We have solved $N=278$
$(N=\sum\limits^{10}_{i=1}\bar N_i, i\neq 4$ and $i\neq 8)$ different
systems of linear algebraic equations at different
orders\footnote[1]{lower index $k$ of order $m_k$ indicates the number of an
example from the set of given-type examples.} $m_k\ (k =1,2,3,4,4',5),$
which is presented below.  {\it A system of the type} $C_2X=Y$. Examples
$1\div5$ from [2,3,16] (see \S 5 Appendix).  $\{m_1:\ 10; 20.\ m_2:\ 3.\
m_3:\ 5; 10; 15.\ m_5:\ 10; 20.$ Here $1<\mu(C^{(m_{\st
k})})\leq1/\sqrt{\varepsilon_1}\},\ \{m_1:\ 30; 40; 50.\ m_2:  4; 5; 6.\
m_3:\ 20; 25; 30; 35.\ m_5:\ 30; 40.$ Here $1/\sqrt{\varepsilon_1}
<\mu(C^{(m_{\st k})})\leq1/\varepsilon_1\},\ \{m_1:\ 60; 70; ...; 100; 150;
200.\ m_2: 8; 9; 10.\ m_3: 40; 45.\ m_4: 5; 6; ...; 18.\ m_5: 50; 60; ...;
100; 150; 200; 300; 400; 500.$ Here $1/\varepsilon_1<\mu(C^{(m_{\st
k})})\}$.  {\it A system of the type} $C_3X=Y$. Examples\footnote[2]{Example
$4'$ is example 4 from [2] (system 9, see \S 5 Appendix),
but with $\varepsilon_0^*$ =0,00000001.}. $6\div10$ from [2,3,16] (see \S 5
Appendix).\newpage
\noindent $\{m_1:$ 10; 20; ...; 100; 150; 200; 300; ...; 900.  $m_2:$ 10;
20; ...; 100; 150; 200; 300; ...; 900. $m_3:$ 10; 21; 3 0; 40; 51; 60; 70;
81; 90; 100; 151; 201; 300; 400; 501; 600; 700; 801; 900. $m_4:$ 10; 30; 40;
60; 70; 90; 100; 300; 400; 600; 700; 900. $m_{4'} :$ 10; 30; 40; 60; 70; 90;
100; 300; 400; 600; 700; 900. $m_5:$ 10; 20; ...; 100. Here
$1<\mu(C^{(m_{\st k})})\leq1/ \sqrt{\varepsilon_1}\},\ \{m_4 :$ 20; 50; 80;
150; 200; 500; 800.  $m_5:$ 150; 200.  Here
$1/\sqrt{\varepsilon_1}<\mu(C^{(m_{\st k})})\leq1/\varepsilon_1\},\ \{m_{4'}
:$ 20; 50; 80; 150; 200; 500; 800. Here $1/\varepsilon_1<\mu(C^{(m_{\st
k})})\}$.  {\it A system of the type} $AX=Y\ (A\ne A^T)$. Examples $11\div15$
from [2,3,16] (see \S 5 Appendix).  $\{m_1:$ 10; 20; ...; 100; 150; 200;
300; 400; 500. $m_2:$ 5. $m_3:$ 10; 20; ...; 100; 150; 200; 250; 300; 400;
500.  $m_4:$ 3; 4; ...; 10. $m_5:$ 5; 6. Here $1<\mu(C^{(m_{\st
k})})\leq1/\sqrt{\varepsilon_1}\},\ \{m_2:$ 6; 7; ...; 10.  $m_4:$ 11; 12;
...; 17. $m_5:$ 7; 8; ...;11.  Here
$1/\sqrt{\varepsilon_1}<\mu(C^{(m_{\st k})})\leq1/\varepsilon_1\},\ \{m_2:$ 11;
12. $m_4:$ 18. $m_5:$ 12. Here $1/\varepsilon_1<\mu(C^{(m_{\st k})})\}$.
{\it A system of the type} $AX=Y\ (A=A^T)$. Examples $16\div20$ from [2,3,16]
(see  \S 5 Appendix).  $\{m_1:$ 10; 20; ...; 100; 150; 200; 300. $m_2:$ 5;
6.  $m_3:$ 5; 6. $m_4:$ 10; 20; ...; 100; 150; 200; 300; 400; 500. Here
$1<\mu(C^{(m_{\st k})}) \leq1/\sqrt{\varepsilon_1}\},\ \{m_2:$ 7; 8; ...; 11.
$m_3:$ 7; 8; ...; 11. Here $1/\sqrt{\varepsilon_1}<\mu(C^{(m_{\st
k})})\leq1/\varepsilon_1\},\ \{m_2:$ 12. $m_3:$ 12. $m_5:$ 5; 10; 15; 20; 25;
30; 35. Here $1/\varepsilon_1<\mu(C^{(m_{\st k})})\}$.

Below, in Tables $1\div3,\ 5\div7$ and $9\div11$, we report the obtained
numerical values of the indicated characteristics
$\delta_L,\ \delta_R,\ \delta_{\tilde X}^*,\ \delta_X,$
of approximate $\tilde X$ and exact  $X$ solutions of systems $WX=Y\ (W:\ C_2; C_3;
A\ne A^T; A=A^T)$ at $1<\mu(W)\leq1/\sqrt{\varepsilon_1}$ --- being
well-posed, $1/\sqrt{\varepsilon_1}<\mu(W)\leq1/\varepsilon_1$ --- being
ill-posed, $1/\varepsilon_1<\mu(W)$ --- being pathologically ill-posed of
these systems, respectively. In Tables 4, 8 and 12, we present averaged
results of Tables $1\div3,\ 5\div7$ and $9\div11$.  Note also that
when $1/\varepsilon_1<\mu(W)$, the subprogram {\tt SVD} stops to
work producing information {\tt INF$=-$1}. The program {\tt TRM}
does not work [19], when $m>100$. Tables 9, 12 do not contain $(*****)$
values of $t$(сек.) above 100 sec.

For an easier apprehension of the calculation results reported in Tables
4, 8 and 12, we plot their ``graphic images-figures".

{\bf Remark 5.} In Tables 4, 8 and 12 (as well as in Figures $1\div9$) ,
we also present the averaged results of computations by subprograms \fbox{\tt
MCS} when $W=(A=A^T)$ and \fbox{\tt OSM} when $W=C_3$, which is to be kept
in mind when analyzing the above Tables and Figures. {\it Explanations} to
some Tables and Figures: on the horizontal axis of Figs. $1\div9$,
in units $10^t$ where values of order $t$ are indicated below the axis,
approximate values of quantities $\delta_L(W),\delta_M(W)$ and $\delta_R(W)$
from Tables 4,8 and 12 are plotted. On the vertical axis of these Figures,
we point out the names of programs through which those values have been
obtained. The horizontal axis of Figs.
$11\div13$ represents relative errors of the r.h.s. of a system with the
Hilbert matrix (see \S 5 of Appendix, example 17), whereas on their
vertical axes, we plot the corresponding relative errors found by various
programs. For names of programs, see the notation. Table 13 contains the
numerical results that are drawn in Figs. $11\div13$. Note that:
{$<\!\delta_Y\!\!>=<\!\frac{||\Delta Y||}{||Y||}\!>$} -- the mean value of
the relative error of perturbation of the r.h.s. of the system;
{$<\!\delta_X\!\!>=<\!\frac{||\Delta X||}{||X||}\!>$} -- exact mean values
corresponding to $<\!\delta_Y\!>.$ Numbers in Table 13 written in line with
a program are average values really obtained by this program for
{$<\!\delta_{\tilde X}\!\!>.$} In Fig.10, the matrix, inverse of the Hilbert
matrix, of order $m=14$ is graphically shown. Along the axis $Z$, the values
of elements of this inverse matrix are indicated. As a result, its
complicated structure is easily visualized. Also, numbers of subspaces
separated by the critical-component method are given in the Figure .
{\small\begin{center}
\begin{tabular}{|c|c|c|c|c|c|c|c|}
\multicolumn{8}{c}{Table 1 (0.173E02$<\mu(W)\leq0.$547E08, 3$\leq m\leq$100,
$\bar\mu(W)=0.$533E07, $\bar{N}_1=117$)}\\ \hline
&PR.&t(sec)&$\delta_L$&$\delta_M$&$\delta_R$&$\delta_{\tilde X}^*$&$\delta_X$\\
\hline
&\tt MCC&0.0004&0.238E-11&0.301E-11&0.470E-10&0.136E01&0.136E01\\
$C_2X=Y$&\tt GS&0.0005&0.238E-11&0.301E-11&0.119E-08&0.136E01&0.136E01\\
$N_1=8$&\tt SVD&0.1568&0.261E-11&0.335E-11&0.723E-08&0.136E01&0.136E01\\
$\bar\mu(C_2)=$&\tt QR&0.0036&0.113E-10&0.175E-10&0.405E-08&0.136E01&0.136E01\\
0.576E07&\tt TRM&0.0366&0.448E-05&0.546E-05&0.575E-05&0.136E01&0.136E01\\
\hline
&\tt MCC&0.0051&0.839E-13&0.868E-13&0.137E-12&0.609E01&0.609E01\\
$C_3X=Y$&\tt GS&0.0028&0.994E-13&0.153E-11&0.289E-10&0.609E01&0.609E01\\
$N_2=54$&\tt QR&6.9793&0.187E-11&0.302E-10&0.294E-08&0.609E01&0.609E01\\
$\bar\mu(C_3)=$&\tt SVD&0.3091&0.609E-12&0.505E-10&0.704E-09&0.609E01&0.609E01\\
$0.106E07$&\tt TRM&5.3795&0.238E-07&0.212E-05&0.212E-05&0.609E01&0.609E01\\
&{\tt OSM}&0.0039&0.561E-01&0.164E00&0.850E04&0.655E01&0.609E01\\
\hline
$AX=Y$&\tt QR&0.1776&0.526E-11&0.101E-09&0.122E-07&0.254E01&0.254E01\\
$A\ne A^T$&\tt MCC&1.0925&0.973E-11&0.117E-09&0.461E-08&0.254E01&0.254E01\\
$N_3=31$&\tt SVD&2.8117&0.973E-11&0.117E-09&0.683E-08&0.254E01&0.254E01\\
$\bar\mu(A)=$&\tt GS&0.1091&0.752E-11&0.144E-09&0.256E-08&0.254E01&0.254E01\\
0.549E07&\tt TRM&2.9615&0.364E-04&0.618E-04&0.639E-04&0.254E01&0.254E01\\
\hline
&\tt MCC&0.9636&0.163E-10&0.951E-10&0.178E-07&0.434E01&0.434E01\\
$AX=Y$&\tt SVD&2.5378&0.163E-10&0.951E-10&0.307E-07&0.434E01&0.434E01\\
$A=A^T$&\tt QR&0.1564&0.129E-10&0.160E-09&0.383E-07&0.434E01&0.434E01\\
$N_4=24$&\tt GS&0.0977&0.373E-10&0.179E-09&0.230E-07&0.434E01&0.434E01\\
$\bar\mu(A)=$&\tt MCS&0.9367&0.153E-09&0.429E-09&0.171E-07&0.434E01&0.434E01\\
0.900E07&\tt TRM&2.5976&0.603E-08&0.328E-07&0.453E-07&0.434E01&0.434E01\\
\hline
\end{tabular}\\ \vspace*{2mm}
\begin{tabular}{|c|c|c|c|c|c|c|c|}
\multicolumn{8}{c}{Table 2 (0.968E08$<\mu(W)\leq$0.399E16, 6$\leq m\leq$80,
$\bar\mu(W)=0.$365E15, $\bar{N}_2=43$)}\\ \hline
&PR.&t(sec)&$\delta_L$&$\delta_M$&$\delta_R$&$\delta_{\tilde X}^*$&$\delta_X$\\
\hline
&\tt MCC&0.0008&0.191E-04&0.353E-04&0.289E-02&0.180E01&0.180E01\\
$C_2X=Y$&\tt GS&0.0009&0.191E-04&0.353E-04&0.517E-02&0.180E01&0.180E01\\
$N_1=13$&\tt SVD&0.7276&0.598E-04&0.110E-03&0.259E00&0.180E01&0.180E01\\
$\bar\mu(C_2)=$&\tt QR&0.0200&0.954E-03&0.221E-02&0.276E00&0.180E01&0.180E01\\
0.351E15&\tt TRM&0.3108&0.208E12&0.208E12&0.210E12&0.131E13&0.180E01\\
\hline
&\tt MCC&0.0067&0.638E-08&0.240E-04&0.277E-02&0.754E01&0.754E01\\
$C_3X=Y$&\tt GS&0.0035&0.672E-06&0.213E-03&0.313E-02&0.754E01&0.754E01\\
$N_2=3$&\tt SVD&8.6764&0.335E-05&0.463E-03&0.453E-01&0.754E01&0.754E01\\
$\bar\mu(C_3)=$&\tt QR&0.4031&0.102E-04&0.962E-03&0.960E-02&0.754E01&0.754E01\\
0.586E15&{\tt OSM}&0.0048&0.468E-01&0.141E00&0.136E11&0.793E01&0.754E01\\
&\tt TRM&6.6831&0.148E14&0.148E14&0.562E14&0.125E15&0.754E01\\
\hline
AX=Y&\tt MCC&1.0966&0.200E-03&0.418E-03&0.239E00&0.254E01&0.254E01\\
$A\ne A^T$&\tt SVD&2.8427&0.200E-03&0.418E-03&0.255E00&0.254E01&0.254E01\\
$N_3=17$&\tt GS&0.1095&0.278E-03&0.514E-03&0.473E00&0.255E01&0.254E01\\
$\bar\mu(A)=$&\tt QR&0.1784&0.338E-03&0.615E-03&0.316E00&0.255E01&0.254E01\\
0.271E15&\tt TRM&2.9711&0.133E12&0.133E12&0.261E12&0.167E12&0.254E01\\
\hline
&\tt MCC&0.9660&0.654E-06&0.161E-04&0.215E00&0.435E01&0.435E01\\
$AX=Y$&\tt SVD&2.5544&0.654E-06&0.161E-04&0.467E00&0.435E01&0.435E01\\
$A=A^T$&\tt MCS&0.9389&0.119E-05&0.548E-04&0.189E00&0.435E01&0.435E01\\
$N_4=10$&\tt GS&0.0979&0.166E-05&0.828E-04&0.452E00&0.435E01&0.435E01\\
$\bar\mu(A)=$&\tt QR&0.1568&0.109E-05&0.976E-04&0.734E00&0.435E01&0.435E01\\
0.253E15&\tt TRM&2.6032&0.435E05&0.435E05&0.460E05&0.543E05&0.435E01\\
\hline
\end{tabular}\\ \vspace*{2mm}
\begin{tabular}{|c|c|c|c|c|c|c|c|}
\multicolumn{8}{c}{Table 3 ($\mu(W)>$0.450E16, 5$\leq m\leq$80, $\bar\mu(W)
>$0.450E16, $\bar{N}_3=46$)}\\ \hline
&PR.&t(sec)&$\delta_L$&$\delta_M$&$\delta_R$&$\delta_{\tilde X}^*$&$\delta_X$\\
\hline
&\tt MCC&0.0010&0.133E00&0.179E00&0.952E02&0.259E01&0.253E01\\
$C_2X=Y$&\tt GS&0.0012&0.482E00&0.103E01&0.924E03&0.378E01&0.253E01\\
$N_1=30$&\tt QR&0.0879&0.274E02&0.276E02&0.211E07&0.182E03&0.253E01\\
$\bar\mu(C_2)>$&\tt TRM&1.6955&0.602E55&0.602E55&0.145E75&0.602E56&0.253E01\\
0.450E16&\tt SVD&\multicolumn{6}{c|}{\tt INF $=-$ 1}\\
\hline
&\tt MCC&0.0069&0.298E00&0.828E00&0.126E02&0.887E01&0.683E01\\
$C_3X=Y$&\tt QR&0.1659&0.724E00&0.127E01&0.398E02&0.129E02&0.683E01\\
$N_2=3$&\tt GS&0.0022&0.174E02&0.184E02&0.940E03&0.142E03&0.683E01\\
$\bar\mu(C_3)>$&{\tt OSM}&0.0034&0.221E02&0.231E02&0.140E18&0.178E03&0.683E01\\
0.450E16&\tt TRM&3.1464&0.777E15&0.777E15&0.611E16&0.430E16&0.683E01\\
&\tt SVD&\multicolumn{6}{c|}{\tt INF $= -$ 1}\\
\hline
$AX=Y$&\tt MCC&0.0179&0.109E-02&0.333E-01&0.106E03&0.856E00&0.856E00\\
$A\ne A^T$&\tt GS&0.0022&0.330E-02&0.681E-01&0.193E03&0.861E00&0.856E00\\
$N_3=4$&\tt QR&0.0040&0.610E-02&0.895E-01&0.209E03&0.861E00&0.856E00\\
$\bar\mu(A)>$&\tt TRM&0.0421&0.166E14&0.166E14&0.278E14&0.208E14&0.856E00\\
0.450E16&\tt SVD&\multicolumn{6}{c|}{\tt INF $=-$ 1}\\
\hline
&\tt MCC&0.0479&0.428E-01&0.365E00&0.139E03&0.994E00&0.989E00\\
$AX=Y$&\tt MCS&0.0439&0.197E00&0.545E00&0.149E03&0.107E01&0.989E00\\
$A=A^T$&\tt GS&0.0051&0.104E05&0.104E05&0.297E22&0.483E04&0.989E00\\
$N_4=9$&\tt QR&0.0090&0.261E12&0.261E12&0.332E25&0.123E12&0.989E00\\
$\bar\mu(A)>$&\tt TRM&0.1165&0.457E16&0.457E16&0.593E18&0.213E16&0.989E00\\
0.450E16&\tt SVD&\multicolumn{6}{c|}{\tt INF $=-$ 1}\\
\hline
\end{tabular}\\ \vspace*{2mm}
\begin{tabular}{|c|c|c|c|c|c|c|c|}
\multicolumn{8}{c}{Table 4 (Mean values of characteristics of Tables $1\div3,
\ \bar{N}_4=206$)}\\ \hline
&PR.&t(sec)&$\delta_L$&$\delta_M$&$\delta_R$&$\delta_{\tilde X}^*$&$\delta_X$\\
\hline
&\tt MCC&0.5154&0.712E-11&0.538E-10&0.561E-08&0.358E01&0.358E01\\
$\bar\mu(W)=$&\tt SVD&3.1214&0.763E-11&0.614E-10&0.119E-07&0.358E01&0.358E01\\
0.533E07&\tt GS&0.0525&0.118E-10&0.819E-10&0.669E-08&0.358E01&0.358E01\\
$\bar{N}_1=117$&\tt QR&0.1617&0.752E-11&0.822E-10&0.138E-07&0.358E01&0.358E01\\
&\fbox{\tt MCS}&0.9367&0.153E-09&0.429E-09&0.171E-07&0.434E01&0.434E01\\
&\tt TRM&2.7438&0.102E-04&0.174E-04&0.180E-04&0.358E01&0.358E01\\
&\fbox{\tt OSM}&0.0039&0.561E-01&0.164E00&0.850E04&0.655E01&0.609E01\\
\hline
&\fbox{\tt MCS}&0.9389&0.119E-05&0.548E-04&0.189E00&0.435E01&0.435E01\\
&\tt MCC&0.5175&0.549E-04&0.123E-03&0.115E00&0.406E01&0.406E01\\
$\bar\mu(W)=$&\tt GS&0.0529&0.749E-04&0.211E-03&0.233E00&0.406E01&0.406E01\\
0.365E15&\tt SVD&3.7003&0.660E-04&0.252E-03&0.257E00&0.406E01&0.406E01\\
$\bar{N}_2=43$&\tt QR&0.1896&0.326E-03&0.971E-03&0.334E00&0.406E01&0.406E01\\
&\fbox{\tt OSM}&0.0048&0.468E-01&0.141E00&0.136E11&0.793E01&0.754E01\\
&\tt TRM&3.1420&0.379E13&0.379E13&0.142E14&0.316E14&0.406E01\\
\hline
&\tt MCC&0.0184&0.119E00&0.351E00&0.882E02&0.333E01&0.280E01\\
$\bar\mu(W)>$&\fbox{\tt MCS}&0.0439&0.197E00&0.545E00&0.149E03&0.107E01&0.989E00\\
0.450E16&{\tt OSM}&0.0034&0.221E02&0.231E02&0.140E18&0.178E03&0.683E01\\
$\bar{N}_3=46$&\tt GS&0.0027&0.260E04&0.260E04&0.742E21&0.124E04&0.280E01\\
&\tt QR&0.0667&0.653E11&0.653E11&0.830E24&0.308E11&0.280E01\\
&\tt TRM&1.2501&0.151E55&0.151E55&0.362E74&0.151E56&0.280E01\\
&\tt SVD&\multicolumn{6}{c|}{\tt INF $=-$ 1}\\
\hline
\end{tabular}
%
\begin{tabular}{|c|c|c|c|c|c|c|c|}
\multicolumn{8}{c}{Table 5 (0.250E03$<\!\mu(W)\!\leq$0.323E07,
150$\leq\!m\!\leq$200,$\bar\mu(W)\!=0.$827E08, $\bar{N}_5\!=16$)}\\ \hline
&PR.&t(sec)&$\delta_L$&$\delta_M$&$\delta_R$&$\delta_{\tilde X}^*$&$\delta_X$\\
\hline
&\tt MCC&0.0133&0.344E-12&0.362E-12&0.724E-12&0.788E01&0.788E01\\
$C_3X=Y$&\tt GS&0.0065&0.345E-12&0.363E-12&0.727E-12&0.788E01&0.788E01\\
$N_2=8$&\tt QR&6.0108&0.622E-12&0.777E-12&0.170E-10&0.788E01&0.788E01\\
$\bar\mu(C_3)=$&\tt SVD&31.137&0.138E-10&0.104E-10&0.274E-10&0.788E01&0.788E01\\
0.202E05&{\tt OSM}&0.0093&0.306E00&0.467E00&0.114E01&0.116E02&0.788E01\\
\hline
$AX=Y$&\tt MCC&22.606&0.331E-13&0.311E-09&0.120E-07&0.719E01&0.719E01\\
$A\ne A^T$&\tt GS&2.2529&0.134E-13&0.374E-09&0.275E-07&0.719E01&0.719E01\\
$N_3=4$&\tt SVD&37.433&0.331E-13&0.436E-09&0.189E-07&0.719E01&0.719E01\\
$\bar\mu(A)=$&\tt QR&3.6669&0.236E-13&0.617E-09&0.306E-07&0.719E01&0.719E01\\
0.123E07&&&&&&&\\
\hline
$AX=Y$&\tt MCC&10.312&0.122E-11&0.155E-09&0.121E-06&0.134E02&0.134E02\\
$A=A^T$&\tt SVD&18.469&0.122E-11&0.155E-09&0.364E-06&0.134E02&0.134E02\\
$N_4=4$&\tt MCS&10.181&0.115E-11&0.215E-09&0.100E-06&0.134E02&0.134E02\\
$\bar\mu(A)=$&\tt QR&1.6768&0.808E-12&0.263E-09&0.371E-06&0.134E02&0.134E02\\
0.123E07&\tt GS&1.0537&0.218E-11&0.336E-09&0.230E-06&0.134E02&0.134E02\\
\hline
\end{tabular}
%
\begin{tabular}{|c|c|c|c|c|c|c|c|}
\multicolumn{8}{c}{Table 6 (0.657E11$<\mu(C_3)\leq$0.594E15, 150$\leq m\leq$
200, $\bar\mu(C_3)=0.$308E15, $\bar{N}_6=3$)}\\ \hline
&PR.&t(sec)&$\delta_L$&$\delta_M$&$\delta_R$&$\delta_{\tilde X}^*$&$\delta_X$\\
\hline
&\tt MCC&0.0591&0.117E-05&0.133E-02&0.207E-01&0.136E02&0.136E02\\
$C_3X=Y$&\tt SVD&43.881&0.488E-05&0.179E-02&0.171E01&0.136E02&0.136E02\\
$N_2=3$&\tt GS&0.0075&0.766E-05&0.277E-02&0.324E-01&0.136E02&0.136E02\\
$\bar\mu(C_3)=$&\tt QR&4.5427&0.114E-04&0.336E-02&0.229E00&0.136E02&0.136E02\\
0.308E15&{\tt OSM}&0.0109&0.135E-02&0.414E-01&0.101E12&0.138E02&0.136E02\\
\hline
\end{tabular}
%
\begin{tabular}{|c|c|c|c|c|c|c|c|}
\multicolumn{8}{c}{Table 7 ($\mu(W)>$0.450E16, 150$\leq m\leq$200,
$\bar\mu(W)>$0.450E16, $\bar{N}_7=5$)}\\ \hline
&PR.&t(sec)&$\delta_L$&$\delta_M$&$\delta_R$&$\delta_{\tilde X}^*$&$\delta_X$\\
\hline
$C_2X=Y$&\tt MCC&0.0049&0.615E-04&0.113E-03&0.000E00&0.719E01&0.719E01\\
$N_1=4$&\tt GS&0.0060&0.615E-04&0.113E-03&0.000E00&0.719E01&0.719E01\\
$\bar\mu(C_2)>$&\tt QR&3.3505&0.401E02&0.408E02&0.128E06&0.653E02&0.719E01\\
0.450E16&\tt SVD&\multicolumn{6}{c|}{\tt INF $=-$ 1}\\
\hline
&\tt MCC&0.0278&0.292E+00&0.819E+00&0.258E02&0.183E02&0.141E02\\
$C_3X=Y$&\tt QR&5.3062&0.551E+00&0.837E+00&0.446E02&0.189E02&0.141E02\\
$N_2=1$&\tt GS&0.0079&0.719E+02&0.728E+02&0.649E04&0.103E04&0.141E02\\
$\bar\mu(C_3)>$&{\tt OSM}&0.0133&0.882E+02&0.892E+02&0.571E18&0.126E04&0.141E02\\
0.450E16&\tt SVD&\multicolumn{6}{c|}{\tt INF $=-$ 1}\\
\hline
\end{tabular}
\begin{tabular}{|c|c|c|c|c|c|c|c|}
\multicolumn{8}{c}{Table 8 (Mean values of characteristics of Tables $5\div7,\
\bar N_8=24$)}\\ \hline
&PR.&t(sec)&$\delta_L$&$\delta_M$&$\delta_R$&$\delta_{\tilde X}^*$&$\delta_X$\\
\hline
&\tt MCC&10.977&0.532E-12&0.155E-09&0.443E-07&0.949E01&0.949E01\\
$\bar\mu(W)=$&\tt SVD&29.013&0.502E-11&0.200E-09&0.128E-06&0.949E01&0.949E01\\
0.827E06&\fbox{\tt MCS}&10.181&0.115E-11&0.215E-09&0.100E-06&0.134E02&0.134E02\\
$\bar{N}_5=16$&\tt GS&1.1044&0.846E-12&0.237E-09&0.848E-07&0.949E01&0.949E01\\
&\tt QR&3.7848&0.498E-12&0.294E-09&0.134E-06&0.949E01&0.949E01\\
&\fbox{\tt OSM}&0.0093&0.306E00&0.467E00&0.114E01&0.116E02&0.788E01\\
\hline
&\tt MCC&0.0591&0.175E-05&0.133E-02&0.207E-01&0.136E02&0.136E02\\
$\bar\mu(C_3)=$&\tt SVD&43.881&0.488E-05&0.179E-02&0.171E01&0.136E02&0.136E02\\
0.308E15&\tt GS&0.0075&0.766E-05&0.277E-02&0.324E-01&0.136E02&0.136E02\\
$\bar{N}_6=3$&\tt QR&4.5427&0.114E-04&0.336E-02&0.229E00&0.136E02&0.136E02\\
&\fbox{\tt OSM}&0.0109&0.135E-02&0.414E-01&0.101E12&0.138E02&0.136E02\\
\hline
\end{tabular}
\newpage
\begin{tabular}{|c|c|c|c|c|c|c|c|}
\hline
&\tt MCC&0.0164&0.146E00&0.410E00&0.129E02&0.127E02&0.106E02\\
$\bar\mu(W)>$&\tt QR&4.3284&0.203E02&0.208E02&0.640E05&0.421E02&0.106E02\\
0.450E16&\tt GS&0.0069&0.360E02&0.364E02&0.325E04&0.519E03&0.106E02\\
$\bar{N}_7=5$&\fbox{\tt OSM}&0.0133&0.882E02&0.892E02&0.571E18&0.126E04&0.141E02\\
&\tt SVD&\multicolumn{6}{c|}{\tt INF $=-$ 1}\\
\hline
\end{tabular}
%
\begin{tabular}{|c|c|c|c|c|c|c|c|}
\multicolumn{8}{c}{Table 9 (0.498E03$\leq\mu(W)\leq$0.318E08, 250$\leq m
\leq$900, $\bar\mu(W)=0.$675E07, $\bar{N}_9=39$)}\\ \hline
&PR.&t(sec)&$\delta_L$&$\delta_M$&$\delta_R$&$\delta_{\tilde X}^*$&$\delta_X$\\
\hline
$C_3X=Y$&\tt MCC&0.0524&0.762E-14&0.403E-13&0.168E-11&0.999E01&0.999E01\\
$\bar\mu(C_3)=$&\tt GS&0.0230&0.718E-14&0.444E-13&0.179E-11&0.999E01&0.999E01\\
0.255E06&\tt OSM&0.0324&0.179E02&0.186E02&0.682E75&0.618E02&0.999E01\\
$N_2=31$&&&&&&&\\
\hline
$AX=Y$&\tt MCC&12.156&0.703E-12&0.611E-08&0.115E-05&0.101E02&0.101E02\\
$A\ne A^T$&\tt GS&1.2256&0.227E-11&0.971E-08&0.179E-05&0.101E02&0.101E02\\
$N_3=4$&&&&&&&\\
$N_4=4$&&&&&&&\\
$A=A^T$&\tt MCC&*****&0.256E-11&0.522E-09&0.182E-05&0.188E02&0.188E02\\
$\bar\mu(A)=$&\tt MCS&*****&0.216E-12&0.227E-08&0.174E-05&0.188E02&0.188E02\\
0.999E07&\tt GS&*****&0.703E-12&0.240E-08&0.197E-05&0.188E02&0.188E02\\
\hline
\end{tabular}
%
\begin{tabular}{|c|c|c|c|c|c|c|c|}
\multicolumn{8}{c}{Table 10 (0.540E15$\leq\mu(C_3)\leq$0.637E15, $m:$
500,800, $\bar\mu(C_3)=$0.589E15, $\bar{N}_{10}=2$)}\\ \hline
&PR.&t(sec)&$\delta_L$&$\delta_M$&$\delta_R$&$\delta_{\tilde X}^*$&$\delta_X$\\
\hline
$C_3X=Y$&\tt MCC&0.0901&0.161E-05&0.174E-02&0.415E-01&0.253E02&0.253E02\\
$\bar\mu(C_3)=$&\tt GS&0.0290&0.294E-04&0.765E-02&0.666E-01&0.253E02&0.253E02\\
0.589E15&\tt OSM&0.0394&0.150E-02&0.915E-01&0.290E75&0.253E02&0.253E02\\
$N_2=2$&&&&&&&\\
\hline
\end{tabular}\\ \vspace*{2mm}
\begin{tabular}{|c|c|c|c|c|c|c|c|}
\multicolumn{8}{c}{Table 11 ($\mu(W)>$0.450E16, 250$\leq m\leq800,\
\bar\mu(W)>$ 0.450E16, $\bar{N}_{11}=7$)}\\ \hline
&PR.&t(sec)&$\delta_L$&$\delta_M$&$\delta_R$&$\delta_{\tilde X}^*$&$\delta_X$\\
\hline
$C_2X=Y$&\tt MCC&0.0088&0.615E-04&0.113E-03&0.000E00&0.106E02&0.106E02\\
$\bar\mu(C_2)>$&\tt GS&0.0115&0.615E-04&0.113E-03&0.000E00&0.106E02&0.106E02\\
0.450E16&&&&&&&\\$N_1=4$&&&&&&&\\
\hline
$C_3X=Y$&\tt MCC&0.0583&0.145E00&0.409E00&0.114E03&0.250E02&0.213E02\\
$\bar\mu(C_3)>$&\tt GS&0.0192&0.167E03&0.168E03&0.190E06&0.394E04&0.213E02\\
0.450E16&{\tt OSM}&0.0300&0.215E03&0.216E03&0.209E75&0.500E04&0.213E02\\
$N_2=3$&&&&&&&\\
\hline
\end{tabular}\\ \vspace*{2mm}
\begin{tabular}{|c|c|c|c|c|c|c|c|}
\multicolumn{8}{c}{Table 12 (Mean values of characteristics of Tables
$9\div11,\ \bar N_{12}=48$)}\\ \hline
&PR.&t(sec)&$\delta_L$&$\delta_M$&$\delta_R$&$\delta_{\tilde X}^*$&$\delta_X$\\
\hline
$\bar\mu(W)=$&\tt MCC&*****&0.109E-11&0.221E-08&0.990E-06&0.130E02&0.130E02\\
0.675E07&\fbox{\tt MCS}&*****&0.216E-12&0.227E-08&0.174E-05&0.188E02&0.188E02\\
$\bar{N}_9=39$&\tt GS&*****&0.993E-12&0.404E-08&0.125E-05&0.130E02&0.130E02\\
&\fbox{\tt OSM}&*****&0.179E02&0.186E02&0.682E75&0.618E02&0.999E01\\
\hline
$\bar\mu(C_3)=$&\tt MCC&0.0901&0.161E-05&0.174E-02&0.415E-01&0.253E02&0.253E02\\
0.589E15&\tt GS&0.0290&0.294E-04&0.765E-02&0.666E-01&0.253E02&0.253E02\\
$\bar{N}_{10}=2$&\fbox{\tt OSM}&0.0394&0.150E-02&0.915E-01&0.290E75&0.253E02&0.253E02\\
\hline
$\bar\mu(W)>$&\tt MCC&0.0335&0.725E-01&0.205E00&0.570E02&0.178E02&0.159E02\\
0.450E16&\tt GS&0.0153&0.835E02&0.840E02&0.950E05&0.198E04&0.159E02\\
$\bar{N}_{11}=7$&\fbox{\tt OSM}&0.0300&0.215E03&0.216E03&0.209E75&0.500E04&0.213E02\\
\hline
\end{tabular}\end{center}}
\begin{figure}{\small
\noindent\begin{tabular}{|c|c|c|c|c|c|}
\multicolumn{6}{r}{Table 13}\\
\hline
\multicolumn{2}{|c|}{$3\leq m\leq6$,}&\multicolumn{2}{c|}{$7\leq m\leq11$,}
&\multicolumn{2}{c|}{$12\leq m\leq13$,}\\
\multicolumn{2}{|c|}{0.524E03$\leq\!\mu(A)\!\leq$0.150E08}&
\multicolumn{2}{c|}{0.475E09$\leq\!\mu(A)\leq\!$0.518E15}&
\multicolumn{2}{c|}{$\mu(A)\!>$0.450E16}\\
\hline
\multicolumn{2}{|c|}{$<\!\delta_X\!>=0.00,\;<\!\delta_Y\!>=0.00$}&
\multicolumn{2}{c|}{$<\!\delta_X\!>=0.00,\;<\!\delta_Y\!>=0.00$}&
\multicolumn{2}{c|}{$<\!\delta_X\!>=0.00,\;<\!\delta_Y\!>=0.00$}\\
\hline
\tt MCS&0.384886772512E-10&\tt MCC&0.64196E-04&\tt MCS&0.182E01\\
\tt QR&0.394340547190E-10&\tt SVD&0.64196E-04&\tt MCC&0.214E01\\
\tt GS&0.432343309209E-10&\tt MCS&0.23960E-03&\tt GS&0.270E01\\
\tt MCC&0.502178164479E-10&\tt GS&0.23988E-03&\tt QR&0.403E01\\
\tt SVD&0.502178164479E-10&\tt QR&0.44822E-03&\tt SVD&\tt INF$=-$1\\
\tt TRM&0.334842097124E-07&\tt TRM&0.58231E05&\tt TRM&0.397E08\\
\hline
\multicolumn{2}{|c|}{$<\!\delta_X\!>=0.10,\;<\!\delta_Y\!>=0.054$}&
\multicolumn{2}{c|}{$<\!\delta_X\!>=0.10,\;<\!\delta_Y\!>=0.047$}&
\multicolumn{2}{c|}{$<\!\delta_X\!>=0.10,\;<\!\delta_Y\!>=0.045$}\\
\hline
\tt MCS&0.999999999782E-01&\tt MCC&0.99945E-01&\tt MCS&0.125E01\\
\tt MCC&0.999999999809E-01&\tt SVD&0.99945E-01&\tt MCC&0.277E01\\
\tt SVD&0.999999999809E-01&\tt QR&0.99982E-01&\tt GS&0.333E01\\
\tt GS&0.999999999948E-01&\tt MCS&0.10005E00&\tt QR&0.465E01\\
\tt QR&0.999999999991E-01&\tt GS&0.10005E00&\tt SVD&\tt INF$=-$1\\
\tt TRM&0.100000344252E00&\tt TRM&0.23802E06&\tt TRM&0.442E10\\
\hline
\multicolumn{2}{|c|}{$<\!\delta_X\!>=0.20,\;<\!\delta_Y\!>=0.107$}&
\multicolumn{2}{c|}{$<\!\delta_X\!>=0.20,\;<\!\delta_Y\!>=0.094$}&
\multicolumn{2}{c|}{$<\!\delta_X\!>=0.20,\;<\!\delta_Y\!>=0.090$}\\
\hline
\tt TRM&0.199999724752E00&\tt QR&0.20003E00&\tt MCC&0.219E01\\
\tt MCC&0.199999999988E00&\tt GS&0.20004E00&\tt MCS&0.323E01\\
\tt SVD&0.199999999988E00&\tt MCS&0.20009E00&\tt QR&0.342E01\\
\tt MCS&0.199999999994E00&\tt MCC&0.20016E00&\tt GS&0.352E01\\
\tt GS&0.199999999997E00&\tt SVD&0.20016E00&\tt SVD&\tt INF$=-$1\\
\tt QR&0.199999999999E00&\tt TRM&0.10208E06&\tt TRM&0.139E11\\
\hline
\multicolumn{2}{|c|}{$<\!\delta_X\!>=0.30,\;<\!\delta_Y\!>=0.161$}&
\multicolumn{2}{c|}{$<\!\delta_X\!>=0.30,\;<\!\delta_Y\!>=0.142$}&
\multicolumn{2}{c|}{$<\!\delta_X\!>=0.30,\;<\!\delta_Y\!>=0.135$}\\
\hline
\tt MCC &0.299999999983E00 &\tt MCS &0.29990E00&\tt MCS&0.172E01\\
\tt SVD &0.299999999983E00 &\tt GS  &0.29991E00&\tt MCC&0.279E01\\
\tt MCS &0.299999999985E00 &\tt QR  &0.29995E00&\tt GS &0.350E01\\
\tt GS  &0.299999999992E00 &\tt MCC &0.29996E00&\tt QR &0.465E01\\
\tt QR  &0.299999999996E00 &\tt SVD &0.29996E00&\tt SVD&\tt INF$=-$1\\
\tt TRM &0.300000398356E00 &\tt TRM &0.60699E04&\tt TRM&0.487E10\\
\hline
\multicolumn{2}{|c|}{$<\!\delta_X\!>=0.39,\;<\!\delta_Y\!>=0.209$}&
\multicolumn{2}{c|}{$<\!\delta_X\!>=0.39,\;<\!\delta_Y\!>=0.184$}&
\multicolumn{2}{c|}{$<\!\delta_X\!>=0.39,\;<\!\delta_Y\!>=0.176$}\\
\hline
\tt MCC&0.389999999988E00&\tt QR &0.38988E00&\tt MCS&0.333E01\\
\tt SVD&0.389999999988E00&\tt MCS&0.39002E00&\tt MCC&0.377E01\\
\tt GS &0.389999999997E00&\tt MCC&0.39003E00&\tt GS &0.392E01\\
\tt MCS&0.389999999998E00&\tt SVD&0.39003E00&\tt QR &0.399E01\\
\tt QR &0.389999999999E00&\tt GS &0.39004E00&\tt SVD&\tt INF$=-$1\\
\tt TRM&0.390000429338E00&\tt TRM&0.31130E06&\tt TRM&0.137E12\\
\hline
\multicolumn{2}{|c|}{$<\!\delta_X\!>=0.60,\;<\!\delta_Y\!>=0.320$}&
\multicolumn{2}{c|}{$<\!\delta_X\!>=0.60,\;<\!\delta_Y\!>=0.282$}&
\multicolumn{2}{c|}{$<\!\delta_X\!>=0.60,\;<\!\delta_Y\!>=0.269$}\\
\hline
\tt MCC&0.597599999976E00&\tt MCC&0.59754E00&\tt MCS&0.197E01\\
\tt SVD&0.597599999976E00&\tt SVD&0.59754E00&\tt MCC&0.334E01\\
\tt MCS&0.597599999983E00&\tt QR &0.59761E00&\tt QR &0.430E01\\
\tt QR &0.597599999985E00&\tt GS &0.59762E00&\tt GS &0.471E01\\
\tt GS &0.597599999990E00&\tt MCS&0.59763E00&\tt SVD&\tt INF$=-$1\\
\tt TRM&0.597600286570E00&\tt TRM&0.19820E06&\tt TRM&0.377E11\\
\hline
\end{tabular}}\end{figure}

%
%

\begin{figure}{\small%
\begin{picture}(50.00,17.00)
\put(0.00,10.00){\framebox(120.00,20.00)[lc]{}}
\put(25.00,25.00){\makebox(0,0)[lc]{The notation used in Figures 1-9:}}
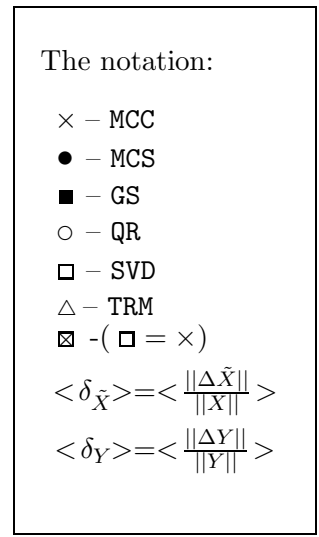
\put(27.0,16.0){\frm{}}
\put(30.0,17.0){\makebox(0,0)[lc]{ -- $\delta^{}_L(W);$}}
\put(53.0,16.0){\fra}
\put(55.0,17.0){\makebox(0,0)[lc]{ -- $\delta^{}_M(W);$}}
\put(77.0,16.0){\frm{$\times$}}
\put(79.0,17.0){\makebox(0,0)[lc]{ -- $\delta^{}_R(W).$}}
\end{picture}
\vspace{-1mm}

\begin{picture}(125.00,45.00)
\multiput(9.00,5.00)(0,5){8}{\line(1,0){1}}
\put(10.00,5.00){\vector(0,1){40.00}}
\multiput(10.00,4.00)(5,0){22}{\line(0,1){1}}
\put(10.00,5.00){\vector(1,0){115.00}}
\put(0.00,0.00){\makebox(0,0)[lc]{$t=$}}
\put(7.00,10.00){\makebox(0,0)[rc]{\tt MCC}}
\put(7.00,15.00){\makebox(0,0)[rc]{\tt SVD}}
\put(7.00,20.00){\makebox(0,0)[rc]{\tt GS}}
\put(7.00,25.00){\makebox(0,0)[rc]{\tt QR}}
\put(7.00,30.00){\makebox(0,0)[rc]{\fbox{\tt MCS}}}
\put(7.00,35.00){\makebox(0,0)[rc]{\tt TRM}}
\put(7.00,40.00){\makebox(0,0)[rc]{\fbox{\tt OSM}}}
\put(14.10,9.10){\frm{}}         \put(24.10,14.10){\frm{}}
\put(34.10,9.10){\fra}           \put(39.10,14.10){\fra}
\put(64.10,9.10){\frm{$\times$}} \put(74.10,14.10){\frm{$\times$}}
\put(29.10,19.10){\frm{}}        \put(19.10,24.10){\frm{}}
\put(44.10,19.10){\fra}          \put(49.10,24.10){\fra}
\put(69.10,19.10){\frm{$\times$}}\put(79.10,24.10){\frm{$\times$}}
\put(54.10,29.10){\frm{}}        \put(89.10,34.10){\frm{}}
\put(59.10,29.10){\fra}          \put(94.10,34.10){\fra}
\put(84.10,29.10){\frm{$\times$}}\put(99.10,34.10){\frm{$\times$}}
\put(104.10,39.10){\frm{}}
\put(109.10,39.10){\fra}
\put(114.10,39.10){\frm{$\times$}}
\put(15.00,10.00){\line(2,1){10.00}}
\put(25.00,15.00){\line(1,1){5.00}}
\put(29.50,20.00){\line(-2,1){9.00}}
\put(20.90,25.00){\line(6,1){34.00}}
\put(54.10,30.00){\line(6,1){36.00}}
\put(91.00,36.00){\line(3,1){14.00}}
{\thicklines
\put(35.00,10.00){\line(1,1){15.00}}
\put(50.00,25.00){\line(2,1){11.00}}
\put(61.00,29.60){\line(6,1){34.00}}
\put(95.00,35.00){\line(3,1){16.00}}}
\put(65.00,10.00){\line(2,1){10.00}}
\put(75.00,15.00){\line(-1,1){5.00}}
\put(70.00,20.00){\line(2,1){10.00}}
\put(80.00,25.00){\line(1,1){5.00}}
\put(85.00,30.00){\line(3,1){15.00}}
\put(100.00,35.00){\line(3,1){15.00}}
\put(10.00,0.00){\makebox(0,0)[cc]{-13}}\put(125.00,2.00){\makebox(0,0)[cc]{$10^t$}}
\multiput(25.00,3.00)(0,2){6}{\line(0,1){1}}   \put(25.00,0.00){\makebox(0,0)[cc]{-12}}
\multiput(50.00,3.00)(0,2){11}{\line(0,1){1}}  \put(50.00,0.00){\makebox(0,0)[cc]{-11}}
\multiput(60.00,3.00)(0,2){14}{\line(0,1){1}}  \put(60.00,0.00){\makebox(0,0)[cc]{-10}}
\multiput(70.00,3.00)(0,2){9}{\line(0,1){1}}   \put(70.00,0.00){\makebox(0,0)[cc]{-9}}
\multiput(85.00,3.00)(0,2){14}{\line(0,1){1}}  \put(85.00,0.00){\makebox(0,0)[cc]{-8}}
\multiput(100.00,3.00)(0,2){16}{\line(0,1){1}} \put(100.00,0.00){\makebox(0,0)[cc]{-5}}
\multiput(105.00,3.00)(0,2){19}{\line(0,1){1}} \put(105.00,0.00){\makebox(0,0)[cc]{-2}}
\multiput(110.00,3.00)(0,2){19}{\line(0,1){1}} \put(110.00,0.00){\makebox(0,0)[cc]{-1}}
\multiput(115.00,3.00)(0,2){19}{\line(0,1){1}} \put(115.00,0.00){\makebox(0,0)[cc]{3}}
\end{picture}

\vspace*{3mm}

\noindent Fig.1. At $\bar\mu(W)\!=$0.533E07 ---
being well-posed, $\bar{N}_1=117$.\\
%
%
\begin{picture}(125.00,48.00)
\multiput(9.00,5.00)(0,5){8}{\line(1,0){1}}
\put(10.00,5.00){\vector(0,1){40.00}}
\multiput(10.00,4.00)(5,0){21}{\line(0,1){1}}
\put(10.00,5.00){\vector(1,0){115.00}}
\put(0.00,0.00){\makebox(0,0)[lc]{$t=$}}
\put(7.00,10.00){\makebox(0,0)[rc]{\fbox{\tt MCS}}}
\put(7.00,15.00){\makebox(0,0)[rc]{\tt MCC}}
\put(7.00,20.00){\makebox(0,0)[rc]{\tt GS}}
\put(7.00,25.00){\makebox(0,0)[rc]{\tt SVD}}
\put(7.00,30.00){\makebox(0,0)[rc]{\tt QR}}
\put(7.00,35.00){\makebox(0,0)[rc]{\fbox{\tt OSM}}}
\put(7.00,40.00){\makebox(0,0)[rc]{\tt TRM}}
\put(14.10,9.10){\frm{}}         \put(24.10,14.10){\frm{}}
\put(19.10,9.10){\fra}           \put(39.10,14.10){\fra}
\put(79.10,9.10){\frm{$\times$}} \put(69.10,14.10){\frm{$\times$}}
\put(34.10,19.10){\frm{}}        \put(29.10,24.10){\frm{}}
\put(44.10,19.10){\fra}          \put(49.10,24.10){\fra}
\put(84.10,19.10){\frm{$\times$}}\put(89.10,24.10){\frm{$\times$}}
\put(54.10,29.10){\frm{}}        \put(64.10,34.10){\frm{}}
\put(59.10,29.10){\fra}          \put(74.10,34.10){\fra}
\put(94.10,29.10){\frm{$\times$}}\put(99.10,34.10){\frm{$\times$}}
\put(104.10,39.10){\fra}
\put(109.10,39.10){\frm{$\times$}}
\put(15.00,10.00){\line(2,1){20.00}}
\put(35.00,20.00){\line(-1,1){5.00}}
\put(30.00,25.00){\line(5,1){26.00}}
\put(56.00,30.33){\line(2,1){9.00}}
\put(65.00,35.00){\line(6,1){30.00}}
\put(95.00,40.00){\line(1,0){10.00}}
{\thicklines
\put(19.00,10.00){\line(4,1){21.00}}
\put(40.00,15.33){\line(1,1){10.00}}
\put(50.00,25.33){\line(2,1){10.00}}
\put(60.00,30.33){\line(3,1){15.00}}
\put(75.00,35.00){\line(6,1){30.00}}}
\put(80.00,10.00){\line(-2,1){10.00}}
\put(70.00,15.00){\line(3,1){15.00}}
\put(85.00,20.00){\line(1,1){15.00}}
\put(100.00,35.00){\line(2,1){10.00}}
\put(10.00,0.00){\makebox(0,0)[cc]{-7}}\put(125.00,2.00){\makebox(0,0)[cc]{$10^t$}}
\multiput(15.00,3.00)(0,2){4}{\line(0,1){1}}  \put(15.00,0.00){\makebox(0,0)[cc]{-6}}
\multiput(35.00,3.00)(0,2){9}{\line(0,1){1}}  \put(35.00,0.00){\makebox(0,0)[cc]{-5}}
\multiput(60.00,3.00)(0,2){14}{\line(0,1){1}} \put(60.00,0.00){\makebox(0,0)[cc]{-4}}
\multiput(65.00,3.00)(0,2){16}{\line(0,1){1}}  \put(65.00,0.00){\makebox(0,0)[cc]{-2}}
\multiput(95.00,3.00)(0,2){14}{\line(0,1){1}} \put(95.00,0.00){\makebox(0,0)[cc]{-1}}
\multiput(100.00,3.00)(0,2){16}{\line(0,1){1}}\put(100.00,0.00){\makebox(0,0)[cc]{10}}
\multiput(105.00,3.00)(0,2){19}{\line(0,1){1}}\put(105.00,0.00){\makebox(0,0)[cc]{12}}
\multiput(110.00,3.00)(0,2){19}{\line(0,1){1}}\put(110.00,0.00){\makebox(0,0)[cc]{13}}
\end{picture}
\vspace*{3mm}

\noindent Fig. 2. At $\bar\mu(W)\!=$0.365E15 ---
being ill-posed, $\bar{N}_2=43$.\\
%
%
\begin{picture}(125.00,43.00)
\multiput(9.00,5.00)(0,5){7}{\line(1,0){1}}
\put(10.00,5.00){\vector(0,1){35.00}}
\multiput(10.00,4.00)(5,0){16}{\line(0,1){1}}
\put(10.00,5.00){\vector(1,0){115.00}}
\put(0.00,0.00){\makebox(0,0)[lc]{$t=$}}
\put(7.00,10.00){\makebox(0,0)[rc]{\tt MCC}}
\put(7.00,15.00){\makebox(0,0)[rc]{\fbox{\tt MCS}}}
\put(7.00,20.00){\makebox(0,0)[rc]{\fbox{\tt OSM}}}
\put(7.00,25.00){\makebox(0,0)[rc]{\tt GS}}
\put(7.00,30.00){\makebox(0,0)[rc]{\tt QR}}
\put(7.00,35.00){\makebox(0,0)[rc]{\tt TRM}}
%
\put(14.10,9.10){\frm{}}         \put(19.10,14.10){\frm{}}
\put(24.10,9.10){\fra}           \put(29.10,14.10){\fra}
\put(44.10,9.10){\frm{$\times$}} \put(49.10,14.10){\frm{$\times$}}
\put(34.10,19.10){\frm{}}
\put(39.10,19.10){\fra}          \put(54.10,24.10){\fra}
\put(64.10,19.10){\frm{$\times$}}\put(69.10,24.10){\frm{$\times$}}
\put(59.10,29.10){\fra}          \put(79.10,34.10){\fra}
\put(74.10,29.10){\frm{$\times$}}\put(84.10,34.10){\frm{$\times$}}
\put(15.00,10.00){\line(1,1){5.00}}
\put(20.00,15.00){\line(3,1){15.00}}
\put(35.00,20.00){\line(4,1){20.00}}
\put(55.00,25.00){\line(1,1){5.00}}
\put(60.00,30.00){\line(4,1){20.00}}
{\thicklines
\put(25.00,10.00){\line(1,1){5.00}}
\put(30.00,15.00){\line(2,1){10.00}}
\put(40.00,20.00){\line(3,1){15.00}}}
\put(45.00,10.00){\line(1,1){5.00}}
\put(50.00,15.00){\line(3,1){15.00}}
\put(65.00,20.00){\line(1,1){10.00}}
\put(75.00,30.00){\line(2,1){10.00}}
\put(10.00,0.00){\makebox(0,0)[cc]{-2}}\put(125.00,2.00){\makebox(0,0)[cc]{$10^t$}}
\multiput(30.00,3.00)(0,2){6}{\line(0,1){1}}  \put(30.00,0.00){\makebox(0,0)[cc]{-1}}
\multiput(45.00,3.00)(0,2){4}{\line(0,1){1}}  \put(45.00,0.00){\makebox(0,0)[cc]{1}}
\multiput(50.00,3.00)(0,2){6}{\line(0,1){1}} \put(50.00,0.00){\makebox(0,0)[cc]{2}}
\multiput(55.00,3.00)(0,2){11}{\line(0,1){1}}  \put(55.00,0.00){\makebox(0,0)[cc]{3}}
\multiput(60.00,3.00)(0,2){14}{\line(0,1){1}} \put(60.00,0.00){\makebox(0,0)[cc]{10}}
\multiput(65.00,3.00)(0,2){9}{\line(0,1){1}}\put(65.00,0.00){\makebox(0,0)[cc]{17}}
\multiput(70.00,3.00)(0,2){11}{\line(0,1){1}}\put(70.00,0.00){\makebox(0,0)[cc]{20}}
\multiput(75.00,3.00)(0,2){14}{\line(0,1){1}}\put(75.00,0.00){\makebox(0,0)[cc]{23}}
\multiput(80.00,3.00)(0,2){16}{\line(0,1){1}}\put(80.00,0.00){\makebox(0,0)[cc]{54}}
\multiput(85.00,3.00)(0,2){16}{\line(0,1){1}}\put(85.00,0.00){\makebox(0,0)[cc]{73}}
\end{picture}

\vspace*{1mm}

\noindent Fig. 3. At $\bar\mu(W)\!>$0.450E16 ---
beind pathologically ill-posed, $\bar{N}_3=46$.\\
%
%
\begin{picture}(125.00,43.00)
\multiput(9.00,5.00)(0,5){7}{\line(1,0){1}}
\put(10.00,5.00){\vector(0,1){35.00}}
\multiput(10.00,4.00)(5,0){19}{\line(0,1){1}}
\put(10.00,5.00){\vector(1,0){115.00}}
\put(0.00,0.00){\makebox(0,0)[lc]{$t=$}}
\put(7.00,10.00){\makebox(0,0)[rc]{\tt MCC}}
\put(7.00,15.00){\makebox(0,0)[rc]{\tt SVD}}
\put(7.00,20.00){\makebox(0,0)[rc]{\fbox{\tt MCS}}}
\put(7.00,25.00){\makebox(0,0)[rc]{\tt GS}}
\put(7.00,30.00){\makebox(0,0)[rc]{\tt QR}}
\put(7.00,35.00){\makebox(0,0)[rc]{\fbox{\tt OSM}}}
\put(19.10,9.10){\frm{}}         \put(34.10,14.10){\frm{}}
\put(39.10,9.10){\fra}           \put(44.10,14.10){\fra}
\put(64.10,9.10){\frm{$\times$}} \put(79.10,14.10){\frm{$\times$}}
\put(29.10,19.10){\frm{}}        \put(24.10,24.10){\frm{}}
\put(49.10,19.10){\fra}          \put(54.10,24.10){\fra}
\put(74.10,19.10){\frm{$\times$}}\put(69.10,24.10){\frm{$\times$}}
\put(14.10,29.10){\frm{}}        \put(89.10,34.10){\frm{}}
\put(59.10,29.10){\fra}          \put(94.10,34.10){\fra}
\put(84.10,29.10){\frm{$\times$}}\put(99.10,34.10){\frm{$\times$}}
\put(20.00,10.00){\line(3,1){15.00}}
\put(35.00,15.00){\line(-1,1){10.00}}
\put(25.00,25.00){\line(-2,1){10.00}}
\put(15.00,30.00){\line(6,1){31.00}}
\put(45.00,35.00){\line(1,0){45.00}}
{\thicklines
\put(40.00,10.00){\line(1,1){19.00}}
\put(59.00,29.00){\line(6,1){36.00}}}
\put(65.00,10.00){\line(3,1){15.00}}
\put(80.00,15.00){\line(-1,1){10.00}}
\put(70.00,25.00){\line(3,1){30.00}}
\put(10.00,0.00){\makebox(0,0)[cc]{-14}}\put(125.00,2.00){\makebox(0,0)[cc]{$10^t$}}
\multiput(25.00,3.00)(0,2){11}{\line(0,1){1}} \put(25.00,0.00){\makebox(0,0)[cc]{-13}}
\multiput(35.00,3.00)(0,2){6}{\line(0,1){1}} \put(35.00,0.00){\makebox(0,0)[cc]{-12}}
\multiput(60.00,3.00)(0,2){14}{\line(0,1){1}} \put(60.00,0.00){\makebox(0,0)[cc]{-10}}
\multiput(70.00,3.00)(0,2){11}{\line(0,1){1}}\put(70.00,0.00){\makebox(0,0)[cc]{-8}}
\multiput(85.00,3.00)(0,2){14}{\line(0,1){1}}\put(85.00,0.00){\makebox(0,0)[cc]{-7}}
\multiput(95.00,3.00)(0,2){16}{\line(0,1){1}} \put(95.00,0.00){\makebox(0,0)[cc]{-1}}
\multiput(100.00,3.00)(0,2){16}{\line(0,1){1}}\put(100.00,0.00){\makebox(0,0)[cc]{0}}
\end{picture}

\vspace*{3mm}

\noindent Fig. 4. At $\bar\mu(W)\!=$0.857E06 ---
being well-posed, $\bar{N}_5=16$.}\end{figure}
%
%
\begin{figure}{\small\noindent\begin{picture}(95.00,35.00)
\multiput(9.00,5.00)(0,5){6}{\line(1,0){1}}
\put(10.00,5.00){\vector(0,1){30.00}}
\multiput(10.00,4.00)(5,0){16}{\line(0,1){1}}
\put(10.00,5.00){\vector(1,0){85.00}}
\put(0.00,0.00){\makebox(0,0)[lc]{$t=$}}
\put(7.00,10.00){\makebox(0,0)[rc]{\tt MCC}}
\put(7.00,15.00){\makebox(0,0)[rc]{\tt SVD}}
\put(7.00,20.00){\makebox(0,0)[rc]{\tt GS}}
\put(7.00,25.00){\makebox(0,0)[rc]{\tt QR}}
\put(7.00,30.00){\makebox(0,0)[rc]{\fbox{\tt OSM}}}
\put(14.10,9.10){\frm{}}         \put(19.10,14.10){\frm{}}
\put(34.10,9.10){\fra}           \put(44.10,14.10){\fra}
\put(59.10,9.10){\frm{$\times$}} \put(74.10,14.10){\frm{$\times$}}
\put(24.10,19.10){\frm{}}        \put(29.10,24.10){\frm{}}
\put(49.10,19.10){\fra}          \put(54.10,24.10){\fra}
\put(64.10,19.10){\frm{$\times$}}\put(79.10,24.10){\frm{$\times$}}
\put(39.10,29.10){\frm{}}
\put(69.10,29.10){\fra}
\put(84.10,29.10){\frm{$\times$}}
\put(15.00,10.00){\line(1,1){15.00}}
\put(30.00,25.00){\line(2,1){10.00}}
{\thicklines
\put(35.00,10.00){\line(2,1){10.00}}
\put(45.00,15.00){\line(1,1){10.00}}
\put(55.00,25.00){\line(3,1){15.00}}}
\put(60.00,10.00){\line(3,1){15.00}}
\put(75.00,15.00){\line(-2,1){10.00}}
\put(65.00,20.00){\line(3,1){15.00}}
\put(80.00,25.00){\line(1,1){5.00}}
\put(10.00,0.00){\makebox(0,0)[cc]{-7}}\put(95.00,2.00){\makebox(0,0)[cc]{$10^t$}}
\multiput(25.00,3.00)(0,2){9}{\line(0,1){1}} \put(25.00,0.00){\makebox(0,0)[cc]{-6}}
\multiput(30.00,3.00)(0,2){11}{\line(0,1){1}}\put(30.00,0.00){\makebox(0,0)[cc]{-5}}
\multiput(55.00,3.00)(0,2){11}{\line(0,1){1}}\put(55.00,0.00){\makebox(0,0)[cc]{-3}}
\multiput(70.00,3.00)(0,2){14}{\line(0,1){1}}\put(70.00,0.00){\makebox(0,0)[cc]{-2}}
\multiput(75.00,3.00)(0,2){6}{\line(0,1){1}}\put(75.00,0.00){\makebox(0,0)[cc]{0}}
\multiput(80.00,3.00)(0,2){11}{\line(0,1){1}}\put(80.00,0.00){\makebox(0,0)[cc]{1}}
\multiput(85.00,3.00)(0,2){14}{\line(0,1){1}}\put(85.00,0.00){\makebox(0,0)[cc]{11}}
\end{picture}\hfill


\vspace*{3mm}

\noindent Fig. 5. At $\bar\mu(W)\!=$0.308E15 ---
being ill-posed, $\bar{N}_6=3$.\\
%
%
\begin{picture}(80.00,33.00)
\multiput(9.00,5.00)(0,5){5}{\line(1,0){1}}
\put(10.00,5.00){\vector(0,1){25.00}}
\multiput(10.00,4.00)(5,0){13}{\line(0,1){1}}
\put(10.00,5.00){\vector(1,0){67.00}}
\put(0.00,0.00){\makebox(0,0)[lc]{$t=$}}
\put(7.00,10.00){\makebox(0,0)[rc]{\tt MCC}}
\put(7.00,15.00){\makebox(0,0)[rc]{\tt QR}}
\put(7.00,20.00){\makebox(0,0)[rc]{\tt GS}}
\put(7.00,25.00){\makebox(0,0)[rc]{\fbox{\tt OSM}}}
%
\put(14.10,9.10){\frm{}}         \put(29.10,14.10){\frm{}}
\put(19.10,9.10){\fra}           \put(34.10,14.10){\fra}
\put(24.10,9.10){\frm{$\times$}} \put(64.10,14.10){\frm{$\times$}}
\put(39.10,19.10){\frm{}}        \put(49.10,24.10){\frm{}}
\put(44.10,19.10){\fra}          \put(54.10,24.10){\fra}
\put(59.10,19.10){\frm{$\times$}}\put(69.10,24.10){\frm{$\times$}}
\put(15.00,10.00){\line(3,1){15.00}}
\put(30.00,15.00){\line(2,1){20.00}}
{\thicklines
\put(20.00,10.00){\line(3,1){15.00}}
\put(35.00,15.00){\line(2,1){20.00}}}
\put(25.00,10.00){\line(1,0){10.00}}
\put(35.00,10.00){\line(6,1){30.00}}
\put(65.00,15.00){\line(-1,1){5.00}}
\put(60.00,20.00){\line(2,1){10.00}}
\put(10.00,0.00){\makebox(0,0)[cc]{-2}}\put(77.00,2.00){\makebox(0,0)[cc]{$10^t$}}
\multiput(20.00,3.00)(0,2){4}{\line(0,1){1}} \put(20.00,0.00){\makebox(0,0)[cc]{-1}}
\multiput(55.00,3.00)(0,2){11}{\line(0,1){1}}\put(55.00,0.00){\makebox(0,0)[cc]{1}}
\multiput(60.00,3.00)(0,2){9}{\line(0,1){1}} \put(60.00,0.00){\makebox(0,0)[cc]{3}}
\multiput(65.00,3.00)(0,2){6}{\line(0,1){1}} \put(65.00,0.00){\makebox(0,0)[cc]{4}}
\multiput(70.00,3.00)(0,2){11}{\line(0,1){1}}\put(70.00,0.00){\makebox(0,0)[cc]{17}}
\end{picture}\hfill
%
%
\begin{picture}(77.00,33.00)
\multiput(9.00,5.00)(0,5){5}{\line(1,0){1}}
\put(10.00,5.00){\vector(0,1){25.00}}
\multiput(10.00,4.00)(5,0){13}{\line(0,1){1}}
\put(10.00,5.00){\vector(1,0){67.00}}
\put(0.00,0.00){\makebox(0,0)[lc]{$t=$}}
\put(7.00,10.00){\makebox(0,0)[rc]{\tt MCC}}
\put(7.00,15.00){\makebox(0,0)[rc]{\fbox{\tt MCS}}}
\put(7.00,20.00){\makebox(0,0)[rc]{\tt GS}}
\put(7.00,25.00){\makebox(0,0)[rc]{\fbox{\tt OSM}}}
\put(24.10,9.10){\frm{}}         \put(14.10,14.10){\frm{}}
\put(29.10,9.10){\fra}           \put(34.10,14.10){\fra}
\put(44.10,9.10){\frm{$\times$}} \put(54.10,14.10){\frm{$\times$}}
\put(19.10,19.10){\frm{}}        \put(59.10,24.10){\frm{}}
\put(39.10,19.10){\fra}          \put(64.10,24.10){\fra}
\put(49.10,19.10){\frm{$\times$}}\put(69.10,24.10){\frm{$\times$}}
\put(25.00,10.00){\line(-2,1){10.00}}
\put(15.00,15.00){\line(1,1){5.00}}
\put(20.00,20.00){\line(6,1){30.00}}
\put(50.00,25.00){\line(1,0){10.00}}
{\thicklines
\put(30.00,10.00){\line(1,1){10.00}}
\put(40.00,20.00){\line(5,1){25.00}}}
\put(45.00,10.00){\line(2,1){10.00}}
\put(55.00,15.00){\line(-1,1){5.00}}
\put(50.00,20.00){\line(4,1){20.00}}
\put(10.00,0.00){\makebox(0,0)[cc]{-14}}\put(77.00,2.00){\makebox(0,0)[cc]{$10^t$}}
\multiput(20.00,3.00)(0,2){9}{\line(0,1){1}}  \put(20.00,0.00){\makebox(0,0)[cc]{-13}}
\multiput(25.00,3.00)(0,2){4}{\line(0,1){1}}  \put(25.00,0.00){\makebox(0,0)[cc]{-12}}
\multiput(40.00,3.00)(0,2){9}{\line(0,1){1}}  \put(40.00,0.00){\makebox(0,0)[cc]{-9}}
\multiput(45.00,3.00)(0,2){4}{\line(0,1){1}}  \put(45.00,0.00){\makebox(0,0)[cc]{-7}}
\multiput(55.00,3.00)(0,2){6}{\line(0,1){1}}  \put(55.00,0.00){\makebox(0,0)[cc]{-6}}
\multiput(65.00,3.00)(0,2){11}{\line(0,1){1}} \put(65.00,0.00){\makebox(0,0)[cc]{1}}
\multiput(70.00,3.00)(0,2){11}{\line(0,1){1}} \put(70.00,0.00){\makebox(0,0)[cc]{74}}
\end{picture}

\vspace*{3mm}

\noindent\begin{tabular}{p{75mm}p{75mm}}
Fig. 6. At  $\bar\mu(W)\!>$0.450E16 --- being pathologically
ill-posed, $\bar{N}_7=5$.&
Fig. 7. At $\bar\mu(W)\!=$0.675E07 ---
being well-posed, $\bar{N}_9=39$.\end{tabular}\\
%
%
\begin{picture}(65.00,28.00)
\multiput(9.00,5.00)(0,5){4}{\line(1,0){1}}
\put(10.00,5.00){\vector(0,1){20.00}}
\multiput(10.00,4.00)(5,0){10}{\line(0,1){1}}
\put(10.00,5.00){\vector(1,0){55.00}}
\put(0.00,0.00){\makebox(0,0)[lc]{$t=$}}
\put(7.00,10.00){\makebox(0,0)[rc]{\tt MCC}}
\put(7.00,15.00){\makebox(0,0)[rc]{\tt GS}}
\put(7.00,20.00){\makebox(0,0)[rc]{\fbox{\tt OSM}}}
\put(14.10,9.10){\frm{}}         \put(19.10,14.10){\frm{}}
\put(29.10,9.10){\fra}           \put(34.10,14.10){\fra}
\put(39.10,9.10){\frm{$\times$}} \put(44.10,14.10){\frm{$\times$}}
\put(24.10,19.10){\frm{}}
\put(49.10,19.10){\fra}
\put(54.10,19.10){\frm{$\times$}}
\put(15.00,10.00){\line(1,1){10.00}}
{\thicklines
\put(30.00,10.00){\line(1,1){5.00}}
\put(35.00,15.00){\line(3,1){15.00}}}
\put(40.00,10.00){\line(1,1){5.00}}
\put(45.00,15.00){\line(2,1){10.00}}
\put(10.00,0.00){\makebox(0,0)[cc]{-7}}\put(65.00,2.00){\makebox(0,0)[cc]{$10^t$}}
\multiput(15.00,3.00)(0,2){4}{\line(0,1){1}} \put(15.00,0.00){\makebox(0,0)[cc]{-6}}
\multiput(20.00,3.00)(0,2){6}{\line(0,1){1}} \put(20.00,0.00){\makebox(0,0)[cc]{-5}}
\multiput(35.00,3.00)(0,2){6}{\line(0,1){1}} \put(35.00,0.00){\makebox(0,0)[cc]{-3}}
\multiput(50.00,3.00)(0,2){9}{\line(0,1){1}} \put(50.00,0.00){\makebox(0,0)[cc]{-2}}
\multiput(55.00,3.00)(0,2){9}{\line(0,1){1}} \put(55.00,0.00){\makebox(0,0)[cc]{74}}
\end{picture}\hspace{18mm}
%
%
\begin{picture}(65.00,28.00)
\multiput(9.00,5.00)(0,5){4}{\line(1,0){1}}
\put(10.00,5.00){\vector(0,1){20.00}}
\multiput(10.00,4.00)(5,0){10}{\line(0,1){1}}
\put(10.00,5.00){\vector(1,0){55.00}}
\put(0.00,0.00){\makebox(0,0)[lc]{$t=$}}
\put(7.00,10.00){\makebox(0,0)[rc]{\tt MCC}}
\put(7.00,15.00){\makebox(0,0)[rc]{\tt GS}}
\put(7.00,20.00){\makebox(0,0)[rc]{\fbox{\tt OSM}}}
\put(14.10,9.10){\frm{}}         \put(29.10,14.10){\frm{}}
\put(19.10,9.10){\fra}           \put(34.10,14.10){\fra}
\put(24.10,9.10){\frm{$\times$}} \put(49.10,14.10){\frm{$\times$}}
\put(39.10,19.10){\frm{}}
\put(44.10,19.10){\fra}
\put(54.10,19.10){\frm{$\times$}}
\put(15.00,10.00){\line(3,1){15.00}}
\put(30.00,15.00){\line(2,1){10.00}}
{\thicklines
\put(20.00,10.00){\line(3,1){15.00}}
\put(35.00,15.00){\line(2,1){10.00}}}
\put(25.00,10.00){\line(5,1){25.00}}
\put(50.00,15.00){\line(1,1){5.00}}
\put(10.00,0.00){\makebox(0,0)[cc]{-3}}\put(65.00,2.00){\makebox(0,0)[cc]{$10^t$}}
\multiput(15.00,3.00)(0,2){4}{\line(0,1){1}} \put(15.00,0.00){\makebox(0,0)[cc]{-2}}
\multiput(20.00,3.00)(0,2){4}{\line(0,1){1}} \put(20.00,0.00){\makebox(0,0)[cc]{-1}}
\multiput(35.00,3.00)(0,2){6}{\line(0,1){1}} \put(35.00,0.00){\makebox(0,0)[cc]{1}}
\multiput(45.00,3.00)(0,2){9}{\line(0,1){1}} \put(45.00,0.00){\makebox(0,0)[cc]{2}}
\multiput(50.00,3.00)(0,2){6}{\line(0,1){1}} \put(50.00,0.00){\makebox(0,0)[cc]{4}}
\multiput(55.00,3.00)(0,2){9}{\line(0,1){1}} \put(55.00,0.00){\makebox(0,0)[cc]{74}}
\end{picture}

\vspace*{3mm}

\noindent%
\begin{tabular}{p{75mm}p{75mm}}
Fig. 8. At $\bar\mu(W)\!=$0.589E15 ---
being ill-posed, $\bar{N}_{10}=2$.&
Fig. 9. At $\bar\mu(W)>$0.450E16 ---
being pathologically ill-posed, $\bar{N}_{11}=7$.\end{tabular}

\vspace*{3mm}


\noindent
\begin{picture}(160,95)
\put(30,115){\special{jpg:graph Gora2.jpg}}
\put(80.00,5.00){\makebox(0,0)[cc]{Fig. 10.}}
\end{picture}}\end{figure}
%
\begin{figure}{\small%
\begin{picture}(70.00,140.00)
\put(7.00,5.00){\vector(0,1){135.00}}\put(7.00,5.00){\vector(1,0){60.00}}
\multiput(4.00,5.00)(0,20){7}{\line(1,0){3}}\multiput(7.00,2.00)(9,0){7}{\line(0,1){3}}
\put(0.00,5.00){\makebox(0,0)[cc]{0.00}}  \put(7.00,0.00){\makebox(0,0)[cc]{0.00}}
\put(0.00,25.00){\makebox(0,0)[cc]{0.10}} \put(16.00,0.00){\makebox(0,0)[cc]{0.05}}
\put(0.00,45.00){\makebox(0,0)[cc]{0.20}} \put(25.00,0.00){\makebox(0,0)[cc]{0.10}}
\put(0.00,65.00){\makebox(0,0)[cc]{0.30}} \put(34.00,0.00){\makebox(0,0)[cc]{0.15}}
\put(0.00,85.00){\makebox(0,0)[cc]{0.40}} \put(44.00,0.00){\makebox(0,0)[cc]{0.20}}
\put(0.00,105.00){\makebox(0,0)[cc]{0.50}}\put(53.00,0.00){\makebox(0,0)[cc]{0.25}}
\put(0.00,125.00){\makebox(0,0)[cc]{0.60}}\put(62.00,0.00){\makebox(0,0)[cc]{0.30}}
\put(0.00,137.00){\makebox(0,0)[cc]{$<\!\delta_{\tilde X}\!\!>$}}
\put(70.00,10.00){\makebox(0,0)[cc]{$<\!\delta_Y\!\!>$}}
\multiput(16.10,5.00)(0,4){34}{\line(0,1){2}}
\multiput(25.50,5.00)(0,4){34}{\line(0,1){2}}
\multiput(35.00,5.00)(0,4){34}{\line(0,1){2}}
\multiput(44.00,5.00)(0,4){34}{\line(0,1){2}}
\multiput(63.00,5.00)(0,4){34}{\line(0,1){2}}
\put(7.00,6.50){\circle*{1.70}} \put(16.10,17.00){\circle*{1.70}}
\put(7.00,8.50){\circle{1.70}}  \put(15.25,18.50){\frm{$\times$}}
\put(6.15,10.00){\fra}          \put(15.25,21.00){\fra}
\put(6.15,12.50){\frm{$\times$}}\put(16.10,24.00){\circle{1.70}}
\put(7.00,15.50){\tri}          \put(16.10,26.00){\tri}
\put(25.50,37.00){\tri}          \put(34.15,56.50){\frm{$\times$}}
\put(24.65,38.50){\frm{$\times$}}\put(35.00,59.70){\circle*{1.70}}
\put(25.50,41.50){\circle*{1.70}}\put(34.15,61.00){\fra}
\put(24.65,43.00){\fra}          \put(35.00,64.50){\circle{1.70}}
\put(25.50,46.00){\circle{1.70}} \put(35.00,67.00){\tri}
\put(43.15,76.50){\frm{$\times$}}\put(62.15,118.00){\frm{$\times$}}
\put(43.15,79.00){\fra}          \put(63.00,121.00){\circle*{1.70}}
\put(44.00,82.00){\circle*{1.70}}\put(63.00,123.00){\circle{1.70}}
\put(44.00,84.00){\circle{1.70}} \put(62.15,124.50){\fra}
\put(44.00,86.00){\tri}          \put(63.00,127.50){\tri}
\end{picture}\hfill
\begin{picture}(70.00,140.00)
\put(7.00,5.00){\vector(0,1){135.00}}\put(7.00,5.00){\vector(1,0){60.00}}
\multiput(4.00,5.00)(0,20){7}{\line(1,0){3}}\multiput(7.00,2.00)(9,0){7}{\line(0,1){3}}
\put(0.00,5.00){\makebox(0,0)[cc]{0.00}}  \put(7.00,0.00){\makebox(0,0)[cc]{0.00}}
\put(0.00,25.00){\makebox(0,0)[cc]{0.10}} \put(16.00,0.00){\makebox(0,0)[cc]{0.05}}
\put(0.00,45.00){\makebox(0,0)[cc]{0.20}} \put(25.00,0.00){\makebox(0,0)[cc]{0.10}}
\put(0.00,65.00){\makebox(0,0)[cc]{0.30}} \put(34.00,0.00){\makebox(0,0)[cc]{0.15}}
\put(0.00,85.00){\makebox(0,0)[cc]{0.40}} \put(44.00,0.00){\makebox(0,0)[cc]{0.20}}
\put(0.00,105.00){\makebox(0,0)[cc]{0.50}}\put(53.00,0.00){\makebox(0,0)[cc]{0.25}}
\put(0.00,125.00){\makebox(0,0)[cc]{0.60}}\put(62.00,0.00){\makebox(0,0)[cc]{0.30}}
\put(0.00,137.00){\makebox(0,0)[cc]{$<\!\delta_{\tilde X}\!\!>$}}
\put(65.00,10.00){\makebox(0,0)[cc]{$<\!\delta_Y\!\!>$}}
\multiput(15.50,5.00)(0,4){34}{\line(0,1){2}}
\multiput(24.50,5.00)(0,4){34}{\line(0,1){2}}
\multiput(33.00,5.00)(0,4){34}{\line(0,1){2}}
\multiput(41.50,5.00)(0,4){34}{\line(0,1){2}}
\multiput(59.00,5.00)(0,4){34}{\line(0,1){2}}
\put(6.15,6.00){\frm{$\times$}} \put(14.65,20.50){\frm{$\times$}}
\put(7.00,9.20){\circle*{1.70}} \put(15.50,23.70){\circle{1.70}}
\put(6.15,10.50){\fra}          \put(15.50,26.00){\circle*{1.70}}
\put(7.00,13.50){\circle{1.70}} \put(14.65,27.50){\fra}
\put(7.00,133.50){\tri}         \put(15.50,139.00){\tri}
\put(24.50,44.00){\circle{1.70}}  \put(33.00,56.70){\circle*{1.70}}
\put(23.65,45.50){\fra}           \put(32.15,58.20){\fra}
\put(24.50,48.50){\circle*{1.70}} \put(33.00,61.40){\circle{1.70}}
\put(23.65,50.00){\frm{$\times$}} \put(32.15,63.00){\frm{$\times$}}
\put(24.50,135.00){\tri}          \put(33.00,133.50){\tri}
\put(41.50,77.00){\circle{1.70}}  \put(58.15,121.00){\frm{$\times$}}
\put(41.50,79.50){\circle*{1.70}} \put(59.00,124.50){\circle{1.70}}
\put(40.65,81.50){\frm{$\times$}} \put(58.15,126.00){\fra}
\put(40.65,84.00){\fra}           \put(59.00,129.50){\circle*{1.70}}
\put(41.50,140.00){\tri}          \put(59.00,134.50){\tri}
\end{picture}\\ \vspace*{1mm}

\noindent
Fig.11 (3$\leq\!m\!\leq$6, 0.524E03$\leq\!\mu(A)\!\leq$0.150E08) \hfill
Fig.12 (7$\leq\!m\!\leq$11, 0.475E09$\leq\!\mu(A)\!\leq$0.518E15)

\vspace*{4mm}

\noindent
\begin{picture}(70.00,70.00)
\put(7.00,5.00){\vector(0,1){65.00}}
\multiput(4.00,5.00)(0,10){6}{\line(1,0){3}}
\put(0.00,5.00){\makebox(0,0)[cc]{0.00}}
\put(0.00,15.00){\makebox(0,0)[cc]{1.00}}
\put(0.00,25.00){\makebox(0,0)[cc]{2.00}}
\put(0.00,35.00){\makebox(0,0)[cc]{3.00}}
\put(0.00,45.00){\makebox(0,0)[cc]{4.00}}
\put(0.00,55.00){\makebox(0,0)[cc]{5.00}}
\put(1.00,63.00){\makebox(0,0)[cc]{$<\!\delta_{\tilde X}\!\!>$}}
\put(7.00,5.00){\vector(1,0){60.00}}
\multiput(7.00,2.00)(9,0){7}{\line(0,1){3}}
\put(7.00,0.00){\makebox(0,0)[cc]{0.00}}
\put(16.00,0.00){\makebox(0,0)[cc]{0.05}}
\put(25.00,0.00){\makebox(0,0)[cc]{0.10}}
\put(34.00,0.00){\makebox(0,0)[cc]{0.15}}
\put(44.00,0.00){\makebox(0,0)[cc]{0.20}}
\put(53.00,0.00){\makebox(0,0)[cc]{0.25}}
\put(62.00,0.00){\makebox(0,0)[cc]{0.30}}
\put(65.00,10.00){\makebox(0,0)[cc]{$<\!\delta_Y\!\!>$}}
\multiput(15.00,5.00)(0,4){16}{\line(0,1){2}}
\multiput(24.00,5.00)(0,4){16}{\line(0,1){2}}
\multiput(32.00,5.00)(0,4){16}{\line(0,1){2}}
\multiput(41.00,5.00)(0,4){16}{\line(0,1){2}}
\multiput(58.00,5.00)(0,4){16}{\line(0,1){2}}
\put(7.00,23.00){\circle*{1.70}} \put(15.00,17.50){\circle*{1.70}}
\put(7.20,27.00){\ti}            \put(15.00,32.70){\ti}
\put(6.15,32.00){\fra}           \put(15.00,38.00){\circle{1.70}}
\put(7.00,45.30){\circle{1.70}}  \put(14.15,51.50){\fra}
\put(7.20,56.50){\tri}           \put(15.00,57.00){\tri}
\put(24.00,27.00){\ti}           \put(32.00,22.00){\circle*{1.70}}
\put(24.00,36.00){\circle*{1.70}}\put(32.00,33.00){\ti}
\put(24.00,38.50){\circle{1.70}} \put(31.15,40.00){\fra}
\put(23.15,40.00){\fra}          \put(32.00,51.50){\circle{1.70}}
\put(24.00,64.00){\tri}          \put(32.00,58.50){\tri}
\put(41.00,38.30){\circle*{1.70}}\put(58.00,25.00){\circle*{1.70}}
\put(41.00,41.00){\ti}           \put(58.00,38.40){\ti}
\put(40.15,42.50){\fra}          \put(58.00,48.00){\circle{1.70}}
\put(41.00,45.50){\circle{1.70}} \put(57.15,52.00){\fra}
\put(41.00,65.00){\tri}          \put(58.00,62.00){\tri}
\end{picture}\hfill
\begin{picture}(40.00,70.00)
\put(15.00,63.00){\makebox(0,0)[cc]{The notation:}}
\put(12.00,55.00){\makebox(0,0)[cc]{$\times\;$-- \tt MCC}}
\put(6.70,50.00){\circle*{1.70}}
\put(6.00,44.00){\fra}
\put(6.70,40.00){\circle{1.70}}
\put(6.00,34.00){\frm{}}
\put(9.50,50.00){\makebox(0,0)[l]{-- \tt MCS}}
\put(9.50,45.00){\makebox(0,0)[l]{-- \tt GS}}
\put(9.50,40.00){\makebox(0,0)[l]{-- \tt QR}}
\put(9.50,35.00){\makebox(0,0)[l]{-- \tt SVD}}
\put(12.00,30.00){\makebox(0,0)[cc]{$\st\triangle\;$-- \tt TRM}}
\put(6.0,25.0){\frm{$\times$}}
\put(10.0,25.0){-( \frm{}}
\put(21.5,26.00){\makebox(0,0)[cc]{$=\times$)}}
\put(20.00,19.00){\makebox(0,0)[cc]{$<\!\delta_{\tilde X}\!\!>=
<\!\frac{||\Delta \tilde X||}{||X||}\!>$}}
\put(20.00,11.00){\makebox(0,0)[cc]{$<\!\delta_Y\!\!>=
<\!\frac{||\Delta Y||}{||Y||}\!>$}}
\put(0.00,0.00){\framebox(40.00,70.00)[cc]{}}
\end{picture}\\ \vspace*{1mm}

\noindent
Fig.13 (12$\leq m\leq$13, $\mu(A)>$0.450E16)}\end{figure}

\newpage

The analysis of numerical results reported in Tables
$1\div13$ and their graphical interpretation with the use of Figs.
$1\div13$ show that our programs
{\tt MCC} and {\tt MCS} provide, on the average, better accuracy
characteristics as compared to the most known analogous programs.

The program {\tt MCC} has also better time characteristics in the case
$W=C_2$, no matter, whether a system of equations is well- or ill- or
pathologocally ill-posed, but {\tt MCC} and {\tt MCS} are about twice as
worse in time as the program {\tt GS} ({\tt DBEQN}) in the case $W=C_3$.
This is owing to the time consumption on the analysis of zeros in computing
$B^{}_{ij}$ --- elements of matrices $B=C^{+}_3$ and on testing various
inequalities in accordance with the algorithm $(2.1)\div(2.5)$. The programs
{\tt MCC} and {\tt MCS} work about 10 times as slow as the program {\tt GS}
({\tt DEQN}) in the general case $W: A=A^T , A\ne A^T$. This is due to
considerable time consumption on reduction of the system $WX=Y$ of the
general form to systems of the type (1.2) and (1.3).

>From the analysis presented it follows that the critical-component method
in its qualitative characteristics is the best one of the methods of
solution of degenerate and ill-posed systems of linear algebraic equations.

\subsubsection{Conclusion}

\abzats
In this paper, we have demonstrated the efficiency of the critical-
component method for numerical solution of degenerate and ill-posed
systems of linear algebraic equations.

We have proved the theorem according to which the only stable normal
solution can surely be obtained for degenerate and ill-posed systems of
linear algebraic equations by the critical-component method.

Results of numerical experiments (278 examples were computed) on the
calculation of basic characteristics of solution of the system $WX=Y$
are presented, and a comparative analysis has been performed, which shows
that the programs {\tt MCC} and {\tt MCS} have, on the average, better
characteristics.

\subsubsection*{References}
\begin{itemize}

\item[{[1]}]Emel'yanenko G.A., Rakhmonov T.T., Dushanov E.B.
Critical-component method for solving systems of linear equations with a
tridiagonal matrix of the general form. JINR preprint, E11-96-105, Dubna, 1996.

\item[{[2]}]Emel'yanenko G.A., Rakhmonov T.T., Dushanov E.B. Algorithms and
programs of the critical-component method of inversion of tridiagonal matrices
and solution of systems of linear equations. JINR preprint, E11-96-106, Dubna, 1996.

\item[{[3]}]Emel'yanenko G.A., Rakhmonov T.T., Dushanov E.B.
Critical-component method solutions of linear algebraic equations. JINR
preprint, E11-96-107, Dubna, 1996.

\item[{[4]}]Tikhonov A.N., Arsenin V.Y. Methods of solution of ill-posed problems.
M.,``Nauka", 1979.

\item[{[5]}]Voevodin V.V. Linear algebra. М., ``Nauka", 1980.

\item[{[6]}]Morozov V.A. Regular methods of solution of  ill-posed problems.
М., ``Nauka", 1987.

\item[{[7]}]Malyshev A.N. Introduction to computational linear algebra
(with application of algorithms on FORTRAN). Novosibirsk,  ``Nauka", 1991.

\item[{[8]}]Voevodin V.V. Computational fundamentals of linear algebra. М.,
``Nauka", 1977.

\item[{[9]}]Voevodin V.V., Kuznetsov U. A. Matrix and computations. М.,
``Nauka", 1984.

\item[{[10]}]Wilkinson J.H., Reinsch C. Reference book of algorithms on ALGOL.
Linear algebra. M., ``Mashinostroenie", 1976.

\item[{[11]}]Seber J. Linear regression analysis. М., ``Mir", 1980.

\item[{[12]}]Emel'yanenko G.A. On properties of systems...(Compact
stable schemes of inversion of tridiagonal mattices). JINR preprint,
Р11-85-488. Dubna, 1985.

\item[{[13]}]Emel'yanenko G.A., Rakhmonov T.T. Quasigeneralized matrix
processes and representations of block-tridiagonal matrices of the
general form, JINR\\ communication, Р11-93-249, Dubna, 1993.

\item[{[14]}]Emel'yanenko G.A. Block-tridiagonal matrices and methods of
numerical solution of spectral problems. Author's summary of doctor's thesis.
VTs so AN USSR. Novosibirsk, JINR, 11-92-4, Dubna, 1992.

\item[{[15]}]Godunov S.K. Solution of systems of linear equations. Novosibirsk,
``Nauka", 1980.

\item[{[16]}]Faddeeva V.N., Kolotilina L.Yu., Computational methods of linear
algebra: set of matrices for testing, parts 1-3 (Materials on mathematical
software of IBM), Leningrad, 1982.

\item[{[17]}]JINR News. Information bulletin of the Joint Institute of Nuclear
Research , No 3, 1996, p.12.

\item[{[18]}]Faddeev D.K., Faddeeva V.N. Computational methods of linear algebra.
М-L., Fizmatgiz, 1963.

\item[{[19]}]Fedorova R.N., Shirokova A.I. Library of programs on FORTRAN.
vv. VI-VII. Description of programs. Dubna, 1990.

\item[{[20]}]Rice J. Matrix computations and mathematical software. M.,
``Mir", 1984.

\item[{[21]}]CERNLIB-CERN Program Library (Short Writeups). Application
Software Group. Computing and Networks Division. CERN, Geneva, Switzerland
(May, 1993).
\end{itemize}

\newpage
\setcounter{subsubsection}{4}
\subsubsection{Appendix}

\abzats
{\it I. Test examples of systems of equations $C_2X=Y$ with two-diagonal
matrices of the general form:}

System 1
$$C_2=\Le 1\ 2\\ \ \quad\cdots\\ \qquad 1\ 2\\ \hspace{3em}1\R,
\quad \begin{array}{l}\x_i=1/i,\\ i=1,2,...,M,\\ \y_i=\frac{3i+1}{
i(i+1)},\ \y_M=1/M,\\ i=1,2,...,M-1;\end{array}\hspace{3.2cm}$$

System 2
$$C_2=\Le\varepsilon^*_0\ r\\ \ \quad\cdots\\ \qquad\varepsilon^*_0\ r\\
\hspace{3em}\varepsilon^*_0\R,\quad\begin{array}{l}\x_i=1/(2i+\varepsilon^*_0),
\ i=1,2,...,M,\\ \y_i=\frac{2i+3\varepsilon^*_0}{(2i+\varepsilon^*_0)(2i+
\varepsilon^*_0+2)},\ \y_M=\varepsilon^*_0/(2M+\varepsilon^*_0),\\ i=1,2,...,M-1,\\
\mbox{ where }r=1-\varepsilon^*_0,\ \varepsilon^*_0=0,01;\end{array}$$

System 3
$$C_2=\Le\frac{7}5\ \frac{11}3\smallskip\\ \ \quad\cdots
\\ \qquad\frac{7}5\ \frac{11}3\smallskip\\ \hspace{3em}\frac{7}5\R,\quad
\begin{array}{l}\x_i=1/(2i+1),\\ i=1,2,...,M,\\ \y_i=\frac{152i+118}
{15(2i+1)(2i+3)},\ \y_M=7/5(2M+1),\\ i=1,2,...,M-1;\end{array}\hspace{5mm}$$

System 4
$$C_2=\Le\varepsilon^*_0\ 2\\ -1\ 2\\ \ \quad\cdots\\ \quad -1\ 2\\ \hspace{2.5em}
\varepsilon^*_1\ 2\\ \hspace{3em}-1\ 2\\ \hspace{4.3em}\cdots\\ \hspace{4em}
-1\ 2\\ \hspace{5.5em}\varepsilon^*_1\R,\quad\begin{array}{l}\x_i=(-1)^{i+1}a,
\\ i=1,2,...,M,\\ \y_1=(\varepsilon^*_0-2)a,\\ \y_i=(-1)^i3a,\\
i=2,3,...,k-1,k+1,...,M-1,\\ \y_k=(-1)^k(2-\varepsilon^*_1)a,\\ \y_M=(-1)^
{M+1}\varepsilon^*_1a,\\ \mbox{ where }a\!=\!1\!+\!\varepsilon^*_0,\
\varepsilon^*_0\!=\!0,0000001,\\ \varepsilon^*_1=0,0001;\end{array}\hspace{1.8cm}$$

System 5
$$C_2=\Le 3\ 7\\  \quad\cdots\\ \qquad 3\ 7\\ \hspace{3em}3\R,
\quad\begin{array}{l}\x_i=1,\\ i=1,2,...,M,\\ \y_i=10,\ \y_M=3,\\
i=1,2,...,M-1.\end{array}\hspace{3.3cm}$$
\newpage

{\it II. Test examples of systems of linear equations $C_3X=Y$ with
tridiagonal matrices of the general form:}

System 6
$$C_3=\Le\ 2 -\!1\\ \!-\!1\ \ 2 -\!1\\ \quad\!\!\ddots\ddots\ddots\\ \quad -1\
\ 2 -\!1\\ \hspace{2.5em}-1\ \ 2\R,\quad\begin{array}{l}\x_i=\frac1i,\
i=1,2,...,M,\ \y_1=2,\ \y_M=\frac{M-2}{M(M-1)},\\ \y_i=\frac{2}{(1-i)i(1+i)},
\ i=2,3,...,M-1;\end{array}$$

System 7
$$C_3=\Le\!-\!1\ \ 1\\ \ 1 -\!2\ \ 1\\ \quad\ddots\ddots\ddots\\ \qquad 1 -\!2\
\ 1\\ \hspace{3.5em}1\ \ a\R,\quad \begin{array}{l}x_i=1+(-1)^i\varepsilon^*_0,
\ i=1,2,...,M,\\ \y_1=2\varepsilon^*_0,\ \y_i=(-1)^{i-1}4\varepsilon^*_0,\
i=2,3,...,M-1,\\ \y_M=(-1)^{M-1}(a+\varepsilon^*_0)\varepsilon^*_0,\\
\mbox{ where }a=\frac{1-M}M,\ \varepsilon^*_0=0,0000001;\end{array}$$

System 8
$$C_3=\Le\ 1 -\!1\\ \!-\!1\ \ 1 -\!1\\ \quad\!\!\ddots\ddots\ddots\\ \quad -1\
\ 1 -\!1\\ \hspace{2.5em}-1\ \ 1\R,\quad \begin{array}{l}x_i=\frac{1}{2i},\
i=1,2,...,M,\ \y_1=\frac14,\ \y_M=\frac{1}{2M(1-M)},\\
\y_i=\frac{i^2+1}{2i(1-i)(1+i)}, i=2,3,...,M-1;\end{array}$$

System 9
$$C_3=\Le 1\quad r\\ p\quad 1\quad r\\ \ \ddots\ddots\ddots\\ \,\ \quad p\quad
1\quad r\\ \hspace{3em}p\quad 1\R,\quad\begin{array}{l}
x_i=1,\ i=1,2,...,M,\ \y_1=2-\varepsilon^*_0,\\ \y_M=2+\varepsilon^*_0,\
\y_i=3,\ i=2,3,...,M-1,\\ \mbox{ where }p=1+\varepsilon^*_0,\
r=1-\varepsilon^*_0,\ \varepsilon^*_0=0,0000001;\end{array}\hspace{4mm}$$

System 10
$$C_3=\Le 6\quad 3\\ 4\quad 6\quad 3\\ \ \ddots\ddots\ddots\\ \,\ \quad 4\quad
6\quad 3\\ \hspace{3em}4\quad 6\R, \begin{array}{l}x_i=1,\ i=1,2,...,M,\\
\y_1=9, \y_M=10,\ \y_i=13,\ i=2,3,...,M-1.\end{array}\hspace{5mm}$$
\newpage

{\it III. Test examples of systems of equations $AX=Y$ with $A\ne A^T$ --
filled matrices of the general form:}

System 11
$$A=\Le\st\;\ M\ \ M-1\ M-2\ \cdots\ 3\ 2\ 333\\ \st M-1\ M-1\ M-2\ \cdots\
3\ 2\ 1\\ \st M-2\ M-2\ M-2\ \cdots\ 3\ 2\ 1\\ \hbox to 3.8cm{{}\leaders
\hbox{$\st\;\ \cdot\;\ $}\hfil {}}\\ \st\;\ 3\quad\;\ 3\quad\;\
3\ \quad\cdots\ 3\ 2\ 1\\ \st\;\ 2\quad\;\ 2\quad\;\ 2\quad\ \cdots\ 2\ 2\ 1
\\ \,\ \varepsilon^*_0\st\quad 1\quad\;\ 1\quad\ \cdots\ 1\ 1\ 1\R,\quad
\begin{array}{l}x_i=1/i,\ i=1,2,...,M,\\ \y_1=\sum\limits^{M-1}_{k=1}\!\!
\frac{M-k+1}k+\frac{333}M,\ \y_M=\sum\limits^M_{k=2}\!\!\frac1k+\varepsilon^*_0,\\
\y_i=(M-i+1)\sum\limits^i_{k=1}\!\!\frac1k+\sum\limits^M_{k=i+1}\!\!
\frac{M-k+1}k,\\i=2,3,...,M-1,\\ \mbox{ where }
\varepsilon^*_0=0,0000001;\end{array}$$

System 12
$$\begin{array}{l}A=(a^{}_{ij}),\ a^{}_{ij}=\frac1{i+j-1},\\i=1,2,...,M-1,\
j=1,2,...,M,\\ a^{}_{M1}=333, \\a^{}_{M1}=\frac1{M+j-1},\ j=2,3,...,M,
\end{array}\quad\begin{array}{l}x_i=1/(2i+1),\ i=1,2,...,M,\\ \y_i=\sum
\limits^M_{k=1}\!\!\frac1{(2k+1)(i+k-1)},\\ i=1,2,...,M-1,\\ \y_M=\sum\limits^
M_{k=2}\!\!\frac1{(2k+1)(i+k-1)}+111;\end{array}$$

System 13
$$\begin{array}{l}A=(a^{}_{ij}),\ a^{}_{1j}=a^{}_{j1}=\frac1{M-j+1},\\
j=1,2,...,M-1,\\ a^{}_{1M}=1+\varepsilon^*_0,\ a^{}_{M1}=1-\varepsilon^*_0,\\a^{}_{ij}
=a^{}_{ji}=\frac1{M-i+1},\\ i=2,3,...,M,\ j=2,3,...,i,\end{array}\quad
\begin{array}{l}x_i=1-\varepsilon^*_0,\ i=1,2,...,M,\\ \y_1=(1-\varepsilon^*_0)
(\sum\limits^{M-1}_{k=1}\!\!\frac1{M-k+1}+1+\varepsilon^*_0),\\ \y_i=
(1-\varepsilon^*_0)(\frac{i}{M-i+1}+\sum\limits^M_{k=1}\!\!\frac1{M-k+1}),\\
i=2,3,...,M-1,\\ \y_M=(1-\varepsilon^*_0)(1-\varepsilon^*_0+\sum\limits^M_{k=2}\!\!
\frac1{M-k+1}),\\ \mbox{ where }\varepsilon^*_0=0,00001;\end{array}$$

System 14
$$A=\Le\st\;\ M\ \ M-1\\ \st M-1\ M-1\ M-2\\ \hbox to 3.5cm{{}\leaders
\hbox{$\st\;\ \cdot\;\ $}\hfil {}}\\ \st\;\ 3\quad\;\ 3\quad\;\ 3\ \quad
\cdots\ 3\ 2\\ \st\;\ 2\quad\;\ 2\quad\;\ 2\quad\ \cdots\ 2\ 2\ 1\\ \st\;\
1\quad\;\ 1\quad\;\ 1\quad\ \cdots\ 1\ 1\ 1\R,\quad\begin{array}{l}\x_i=(-1)^
i/i,\ i=1,2,...,M,\\ \y_i=(1+\frac{(-1)^i(i-M)}{i+1})\sum\limits^i_{k=1}
\frac{(-1)^k}k,\\i=1,2,...,M-1,\\ \y_M=\sum\limits^M_{k=1}(-1)^k/k;
\end{array}\hspace{1cm}$$

System 15
$$\begin{array}{l}A=(a^{}_{ij}),\ a^{}_{ij}=\frac1{i-j+M},\\i=1,2,...,M,\
j=1,2,...,M,\end{array}\quad\begin{array}{l}\x_i=1/i,\ i=1,2,...,M,\\
\y_i=\sum\limits^M_{k=1}\!\!\frac1{k(i-k+M)},\\ i=1,2,...,M.\end{array}$$
\newpage

{\it IV. Test examples of systems of equations $AX=Y$ with $A=A^T$ --
filled matrices of the general form:}

System 16
$$\begin{array}{l}A=(a^{}_{ij}),\ a^{}_{1j}=a^{}_{j1}=\frac1{M-j+1},\\
j=1,2,...,M-1,\\ a^{}_{M1}=a^{}_{1M}=M+\varepsilon^*_0,\\
a^{}_{ij}=a^{}_{ji}=\frac1{M-i+1},\\ i=2,3,...,M,\ j=2,3,...,i,\\ \x_i=(i+1)/i,
\ i=1,2,...,M,\end{array}\quad\begin{array}{l}\y_1=\sum\limits^{M-1}_{k=1}
\frac{k+1}{k(M-k+1)}+\frac{(M+\varepsilon^*_0)(M+1)}M,\\ \y_i=\frac1{M-i+1}\sum
\limits^i_{k=1}\frac{k+1}k+\sum\limits^M_{k=i+1}\frac{k+1}{k(M-k+1)},\\
\y_M=\sum\limits^M_{k=2}\frac{k+1}k+2(M+\varepsilon^*_0),\\
i=1,2,...,M-1,\mbox{ where }\varepsilon^*_0=0,0000001;\end{array}$$

System 17($A$ -- the Hilbert matrix)
$$\begin{array}{l}A=(a^{}_{ij}),\ a^{}_{ij}=\frac1{i+j-1},\\ i=1,2,...,M,\
j=1,2,...,M,\end{array}\quad\begin{array}{l}\x_i=1/i,\ i=1,2,...,M,\\ \y_
i=\sum\limits^M_{k=1}\frac1{k(i+k-1)},\ i=1,2,...,M;\end{array}$$

System 18
$$\begin{array}{l}A=(a^{}_{ij}),\ a^{}_{1j}=a^{}_{j1}=\frac1{i+j-1},\\
j=1,2,...,M-1,\\ a^{}_{M1}\!=\!a^{}_{1M}\!=\!333,\ a^{}_{ij}\!=\!a^{}_{ji}\!=
\!\frac1{i+j-1},\\ i=2,3,...,M,\ j=2,3,...,i,\\ \x_i=(-1)^i/i,\ i=1,2,...,M,
\end{array}\quad\begin{array}{l}\y_1=\sum\limits^{M-1}_{k=1}\frac{(-1)^k}
{k^2}+\frac{(-1)^M333}{M},\\ \y_{i-1}\!=\!\sum\limits^M_{k=1}\!\!\frac{(-1)^k}
{k(k+i-2)},\ i=3,4,...,M,\\ \y_M=\sum\limits^M_{k=2}\frac{(-1)^k}{k(k+M-1)}+
333;\end{array}$$

System 19
$$A=\Le\st\ \ \varepsilon^*_0\ \ M-1\ M-2\ \cdots\ 3\ 2\ M\\ \st M-1\ M-1\ M-2\
\cdots\ 3\ 2\ 1\\ \hbox to 3.8cm{{}\leaders\hbox{$\st\;\ \cdot\;\ $}\hfil {}}\\
\st\;\ 2\quad\;\ 2\quad\;\ 2\quad\ \cdots\ 2\ 2\ 1\\ \st\ M\quad\ 1\quad\;\
1\quad\ \cdots\ 1\ 1\ 1\R,\quad\begin{array}{l}\x_i=1-\varepsilon^*_0,\ i=1,2,...,M,
\\ \y_1=\frac{(1-\varepsilon^*_0)(M^2-M+2\varepsilon^*_0)}2,\\ \y_i=\frac{(1-
\varepsilon^*_0)(M-i)(i+M-3)}2,\\ \y_M=(2M-1)(1-\varepsilon^*_0),\\
i=2,3,...,M-1,\mbox{ where }\varepsilon^*_0=0,0000001;\end{array}$$

System 20(det(A)=0)
$$A=\left[\begin{array}{cccccc}a&b&b&\cdots&b&a\\b&a&b&\cdots&b&b\\ \cdot&
\cdot&\cdot&\cdots&\cdot&\cdot\\b&b&b&\cdots&a&b\\a&b&b&\cdots&b&a\R,\quad
\begin{array}{l}\x_i=\frac{(-1)^i}{2i+1},\ i=1,2,...,M,\\ \y_1\!=\!y_M=b\!\sum
\limits^{M-1}_{k=2}\!\!\frac{(-1)^k}{2k+1}+a\!\left(\!\frac{(-1)^M}{2M+1}-
\frac13\!\right),\\ \y_i=b\sum\limits^M_{k=1}\frac{(-1)^k}{2k+1}-\frac{(-1)^
ia}{2i+1},\\i=2,3,...,M-1,\mbox{ where } a=1-\varepsilon^*_0,\\
b=1+\varepsilon^*_0,\ \varepsilon^*_0=1\cdot 10^{-11}.\end{array}$$
\end{document}